\documentclass[reqno, a4paper]{amsart}
\usepackage{amssymb, mathdots}
\usepackage{mathabx}
\usepackage{mathrsfs}
\usepackage[utf8]{inputenc}
\usepackage{amsmath,  graphicx, tensor, stackrel, mathtools}
\usepackage{enumitem}
\usepackage[all]{xy}
\usepackage{tikz}
\usetikzlibrary{matrix, arrows.meta}
\usetikzlibrary{cd}

\usepackage{hyperref}
\DeclareSymbolFont{bbold}{U}{bbold}{m}{n}
\DeclareSymbolFontAlphabet{\mathbbold}{bbold}
\def\qmod#1#2{{\hbox{}^{\displaystyle{#1}}}\!\big/\!\hbox{}_{
\displaystyle{#2}}}

 \def\psp#1#2%
  {\mathop{}%
   \mathopen{\vphantom{#2}}^{#1}%
   \kern-\scriptspace%
    \hskip -0.3mm{#2}} 
 \def\psb#1#2%
  {\mathop{}%
   \mathopen{\vphantom{#2}}_{#1}%
   \kern-\scriptspace%
    \hskip -0.0mm{#2}} 
 \def\pscr#1#2#3%
  {\mathop{}%
   \mathopen{\vphantom{#3}}^{#1}_{#2}%
   \kern-\scriptspace%
    \hskip -0.3mm{#3}} 

\def\C{{\mathbb C}}

\def\N{{\mathbb N}}

\def\P{{\mathbb P}}

\def\Z{{\mathbb Z}}

\def\map{\longrightarrow}
\def\textmap#1{\mathop{\vbox{\ialign{
                                  ##\crcr
      ${\scriptstyle\hfil\;\;#1\;\;\hfil}$\crcr
      \noalign{\kern 1pt\nointerlineskip}
      \rightarrowfill\crcr}}\;}}
\def\bigtextmap#1{\mathop{\vbox{\ialign{
                                  ##\crcr
      ${\hfil\;\;#1\;\;\hfil}$\crcr
      \noalign{\kern 1pt\nointerlineskip}
      \rightarrowfill\crcr}}\;}}
      
\newcommand{\cal}{\mathcal}
\def\textlmap#1{\mathop{\vbox{\ialign{
                                  ##\crcr
      ${\scriptstyle\hfil\;\;#1\;\;\hfil}$\crcr
      \noalign{\kern-1pt\nointerlineskip}
      \leftarrowfill\crcr}}\;}}
 
\def\ag{{\mathfrak a}}

\def\jg{{\mathfrak j}}

\def\mg{{\mathfrak m}}

\def\pg{{\mathfrak p}}

\def\Ag{{\mathfrak A}}

\def\Rg{{\mathfrak R}}

\def\Tg{{\mathfrak T}}

\newtheorem{sz}{Satz}[section]
\theoremstyle{remark}
\newtheorem{re}[sz]{Remark} 
\theoremstyle{plain}
\newtheorem{thr}[sz]{Theorem}
\newtheorem{pr}[sz]{Proposition}

\newtheorem{dt}[sz]{Definition}
\newtheorem{lm}[sz]{Lemma}

\theoremstyle{definition}

\newtheorem*{main-assumpt}{Main Assumption}

\def\End{\mathrm {End}}
\def\Aut{\mathrm {Aut}}

\def\GL{\mathrm {GL}}

\def\Pic{\mathrm {Pic}}

\def\Hom{\mathrm{Hom}}

\def\id{ \mathrm{id}}
\def\im{\mathrm{im}}

\def\coker{\mathrm{coker}}

\def\Ext{\mathrm{Ext}}

\newcommand\smvee{{\hskip -0.1ex \raise 0.2ex\hbox{$\scriptscriptstyle\vee$}}\hskip -0,3ex}

\def\edf{\coloneq}

\def\Torss{\mathcal{T}{\rm ors}}

\def\Aff{\mathrm{Aff}}
\def\Afff{\mathcal{A}{\rm ff}}
\def\Afffs{\tensor[_S]{\Afff}{}}
\def\Torsss{\tensor[_S]{\Torss}{}}

\begin{document}

\title[The moduli stack of affine bundles]{The moduli stack of affine bundles with fixed linearisation type }

\author{Naoufal Bouchareb \and Andrei Teleman}
\address{Aix Marseille Univ, CNRS, I2M, Marseille, France.}

\email[Naoufal Bouchareb]{noufelbouchareb@hotmail.fr}
\email[Andrei Teleman]{andrei.teleman@univ-amu.fr}
 
\begin{abstract}
The goal of this article is the classification of affine bundles with fixed linearisation type over a compact complex manifold $X$.  More precisely, let $E$ be a fixed holomorphic vector bundle over $X$ and ${\cal E}$ the associated locally free sheaf. We will assume that there exists an effective divisor $D$ on $X$ such that $H^1(X,{\cal E}(D))=0$. This condition is always satisfied when $X$ is a projective complex manifold.

We  define in a natural way the  moduli stack  of affine bundles on $X$ whose linearisation is isomorphic to $E$, and we  prove that this moduli stack is isomorphic to  the quotient stack $\big[H^0(X,{\cal E}(D)_D)/H^0(X,{\cal E}(D))\ltimes \Aut(E)\big]$ of the complex vector space $H^0(X,{\cal E}(D)_D)$ by the semi-direct product of the automorphism group $\Aut(E)$ (which is always an affine algebraic group) by the complex vector space $H^0(X,{\cal E}(D))$. Our isomorphism is obtained by constructing {\it explicitly} comparison functors  $\mathscr{A}$, $\mathscr{T}$ between the two stacks and {\it
explicit} natural isomorphisms  $\mathscr{A}\circ\mathscr{T}\textmap{\simeq} \id$, $\id \textmap{\simeq}\mathscr{T}\circ\mathscr{A}$.

  \end{abstract}


\maketitle

\tableofcontents

\section{Introduction}

\subsection{The main result}

Our goal is the classification of affine bundles with fixed linearisation type over a compact complex manifold $X$ in full generality. The first problem in this direction is the set theoretical classification; this problem is explicitly solved by Theorem \ref{ClassTh} (see also \cite[Theorem 1.9]{Bou1}). This general classification  theorem has been effectively used  in \cite{Bou1} in the special case $X=\P^1_\C$, showing  that, even on the simplest Riemann surface, the classification of affine bundles is much more complicated and interesting than the classification of vector bundles.

It is easy to see that, in general, the set of isomorphism classes of affine bundles with fixed linearisation type does not have a natural complex space structure. In other words we do not have a moduli space in the classical sense parameterising this set.  The correct framework for our classification problem is Stack Theory. We refer to \cite{He}, \cite{Wa} for the general Moduli Stack Theory for principal bundles in the algebraic framework.

Let $X$ be a compact complex manifold, $E$ a holomorphic vector bundle on $X$ and ${\cal E}$ the associated locally free sheaf. We define a holomorphic stack (see \cite{Mee}) denoted $\mathscr{A}{\rm ff}_E$ on the category of complex spaces  in the following way: we associate to a complex space $S$ the groupoid $\Afffs_E$ whose objects are affine $A$ bundles on $S\times X$, whose linearisation $L(A)$ is ``fibrewise isomorphic to $E$" in the following sense: for any point $s\in S$ there exists an open neighbourhood $U$ of $s$ in $S$ and a bundle isomorphism $p_X^{U\times X*}(E)\textmap{\simeq} L(A)|_{U\times X}$ (see section \ref{fbrwse-isom-E-sect} for details). This condition commutes obviously with pull-backs. 

If $S$ is reduced, this condition reduces to 
$$
\forall s\in S,\ L(A)|_{\{s\}\times X}\simeq E 
$$
(via the obvious identification $\{s\}\times X\textmap{\simeq}X$). Using Proposition \ref{AffBdls-AugmBdls} we will identify $\Afffs_E$ with the groupoid of augmented bundles $\tilde {\cal E}\textmap{\varphi}{\cal O}_{S\times X}$ over $S\times X$ with ${\cal E}_\varphi\edf \ker(\varphi)$ fibrewise isomorphic to ${\cal E}$.
 
We will assume  that there exists an effective divisor $D$ on $X$ such that 
$$H^1(X,{\cal E}(D))=0,$$
 we  fix such a divisor and put $V_D\edf H^0(X,{\cal E}(D)_D)$, $W_D\edf H^0(X,{\cal E}(D))$.  The automorphism group $G\edf\Aut(E)$ is always an affine algebraic group (see Proposition \ref{Aut(E)AffAlg}) and acts on the vector spaces $V_D$, $W_D$ in a natural way.  The semidirect product $G_D\edf W_D\ltimes G$ comes with a natural action $\alpha_D$ on $V_D$ induced by the  mentioned above $G$-action on $V_D$  and the action of $W_D$ on $V_D$ by translations via the natural linear map $W_D\to V_D$.

Our main result is
\theoremstyle{plain}
\newtheorem*{equiv-stacks-th-intro}{Theorem \ref{equiv-stacks-th}}
\begin{equiv-stacks-th-intro} 
Let $X$ be a compact complex manifold, and $E$ be a fixed holomorphic vector bundle on $X$ and ${\cal E}$ the associated locally free sheaf. Let $D$ be an effective divisor on $X$ such that $H^1(X,{\cal E}(D))=0$. The  holomorphic stack $\mathscr{A}{\rm ff}_E$ is equivalent to the quotient stack $[V_D/G_D]$.
\end{equiv-stacks-th-intro}

Therefore the stack $\Afff_E$ is equivalent to the quotient stack $[V_D/G_D]$ whose fiber over a complex space  $S$ is the groupoid $\Torsss^{\alpha_D}_{G_D}$ whose objects are pairs $(\Theta,\rho)$ consisting of an ${\cal O}_S(G_D)$-torsor  on $S$ and a ${\cal O}_S(G_D)$-equivariant sheaf morphism $\rho:\Theta\to {\cal O}_S(V_D)$.

Our proof will be obtained by constructing explicitly functors 
$$\mathscr{A}:[V_D/G_D]\to \mathscr{A}{\rm ff}_E,\ \mathscr{T}: \mathscr{A}{\rm ff}_E\to [V_D/G_D]$$
and explicit natural isomorphisms $\mathscr{A}\circ \mathscr{T}\textmap{\simeq} \id_{\mathscr{A}{\rm ff}_E}$, $\mathscr{T}\circ \mathscr{A}\textmap{\simeq} \id_{[V_D/G_D]}$.   The functors $\mathscr{A}$, $\mathscr{T}$  will be defined by constructing, for any individual complex space $S$, functors $\tensor[_S]{\Ag}{}:\Torsss^{\alpha_D}_{G_D}\to \Afffs_E$, $\tensor[_S]{\Tg}{}:\Afffs_E\to \Torsss^{\alpha_D}_{G_D}$  with clear geometric interpretations. The construction of these functors will be explained in the next section. Using the terminology of \cite[p. 11]{Gom}, our result shows that the stack $\mathscr{A}{\rm ff}_E$ is isomorphic to the quotient stack $[V_D/G_D]$.

 The results of this article have been presented in the second part of the first author's PhD thesis \cite{Bou2}. 
 
\subsection{Previous work on holomorphic affine bundles}

In his famous work \cite{We}, Weil introduces and studies affine line bundles over an algebraic curve $C$. Using the terminology of this article, the objects studied by Weil are algebraic affine line bundles over    $C$ whose linearisation is isomorphic to the line bundle associated with a divisor $D$ on $C$.  His classification result \cite[section 4]{We} can be recovered as a special case  of (the algebraic-geometric analogue of) the general set theoretical classification Theorem \ref{ClassTh}  which will be stated in section \ref{section-AffBundles}.  
 
\vspace{2mm}

Affine line bundles in the holomorphic category have been considered in Enoki's article about class VII surfaces with curves \cite{En}. In this article the author proves that any minimal class VII surface with $b_2>0$ which contains a non-zero homologically trivial  divisor can be obtained as a compactification of the total space of an  affine line bundle with negative degree over an elliptic curve. The compactifying divisor is a homologically trivial cycle of rational curves.
\vspace{2mm}

The classification problem of affine bundles in the complex geometric category has been previously studied by Gómez-Presas   and Plechinger. We summarise their main results below:

In \cite{GoPr} the authors state the equivalence between the groupoid of affine bundles of rank $r$ and the groupoid of augmented bundles of rank $r+1$ (i.e.  of epimorphisms $\tilde {\cal E}\textmap{\varphi}{\cal O}_X$, where $\tilde {\cal E}$ is a locally free  sheaf of rank $r+1$) which we will recall briefly below. The authors use this result to construct a bijection between the set of conjugacy classes of semi-simple representations of the fundamental group $\pi_1(X,x_0)$ of a compact Kähler manifold into the affine group $\Aff(r,\C)$, and the set of isomorphism classes of ``Higgs extensions" with polystable Higgs kernel (see \cite[Definition 3.4, Theorem 3.5]{GoPr}).

Plechinger's main results \cite{Pl1}, \cite{Pl2} are dedicated to the classification of affine line bundles on a compact complex manifold:
\begin{itemize} 
\item He obtained classification results for framed affine line bundles on compact complex manifolds. In Plechniger's formalism, an $x_0$-framed affine line bundle on $X$ is a pair $(A,\varepsilon)$, where $A$ is an affine   on $X$ and $\varepsilon$ is  a linear isomorphism between  $\C$ and the fibre $L(A)_{x_0}$ of the linearisation $L(A)$ of $A$ at a fixed point $x_0\in X$. He gives a sufficient and necessary condition for the existence of a moduli space in the sense of classical complex geometry  which classifies $x_0$-framed affine line bundles $A$ on $X$  with fixed Chern class $c_1(L(A))=c\in H^2(X,\Z)$,   up to isomorphism.    If this existence condition is satisfied,  the obtained moduli space is a linear space over the component $\Pic^c(X)$ of the Picard group $\Pic(X)$.  The moduli space of non-framed affine line bundles with $c_1(L(A))=c$ can be recovered as the (highly non-Hausdorff) $\C^*$-quotient of the moduli space of framed line bundles.

\item For the general case, Plechinger constructed and studied  the holomorphic stack which classifies framed  and non-framed holomorphic affine line bundles  which replaces (in the framed case) the classical moduli space when the existence condition is not satisfied. 
  
\end{itemize}

\section{Preliminaries}\label{prelim}

\subsection{Affine bundles}\label{section-AffBundles}

Let $\Aff(r,\C)$ be the group of affine transformations of $\C^r$. Our main objects are holomorphic bundles in the following sense: 
\begin{dt} 
Let $X$ be a complex manifold. A holomorphic affine bundle of rank $r$ on $X$ is a couple $(\pi: A \to X, \mathcal{T})$ consisting of a surjective holomorphic submersion $\pi$ and a bundle  atlas ${\cal T}$ (see for instance \cite{Sch})  with structure group $\Aff(r,\C)$ on $\pi$ which  is maximal with respect to inclusion.	
\end{dt}
The classical theory of  principal bundles in the  differential geometric framework (see for instance \cite{KN}, \cite{Te}) extends naturally in the holomorphic category. 
With any affine bundle $(\pi : A \to  X, \mathcal{T})$ one can associate an $\Aff(r,\C)$-principal  bundle as follows: for any $x\in X$  we define 
\begin{equation*}
\begin{split}
P_{A,x}\coloneq \{f: \C^r \to  A_x|\ \hbox{ for any }\pi^{-1}(U_\tau)\textmap{\tau} U_\tau\times\C^r\in & {\cal T} \hbox{ with }x\in U_\tau, \\ &\hbox{ we have }   \tau_x \circ f \in \Aff(r,\C)\}.
\end{split}
\end{equation*}
The disjoint union $P_A\coloneq \coprod_{x\in X}P_{A,x}$
has a natural complex manifold structure and comes with an obvious $\Aff(r,\C)$-right action such that the natural map $P_A\textmap{\pi_A} X$ becomes a holomorphic $\Aff(r,\C)$-principal bundle.

Conversely, for any holomorphic  principal $\Aff(r,\C)$-bundle $P$ on $X$ we define the  affine bundle  
$$
A^P\edf P\times_{\Aff(r,\C)}\C^r
$$
to be the bundle associated with $P$ and the standard  action of its structure group $\Aff(r,\C)$ on $\C^r$. Similarly to the case of vector bundles, we have 
%
%
\begin{re}
The functors $A\mapsto P_A$, $P\mapsto A^P$ define an   equivalence between the groupoid of affine rank $r$-bundles on $X$   and the groupoid of principal $\Aff(r,\C)$-bundles on $X$.	
\end{re}

The set ${\cal T}$  of local trivialisations of $A$ can be identified with the set of local trivialisations (or, equivalently, of of local sections) of $P_A$.

Let $(\pi:A\to X,{\cal T})$ be an affine bundle on $X$ and $x\in X$. A biholomorphic automorphism $u$ of $A_x$ will be called {\it translation} if there exists (or, equivalently, if for any) $\tau\in {\cal T}$ with $x\in U_\tau$ the composition $\tau_x \circ u \circ \tau_x^{-1}$ is a translation of $\C^r$. The set $L(A)_x$ of translations of $A_x$ has the natural structure of an $r$-dimensional complex vector space, and the natural  action of $L(A)_x$ on $A_x$ is free and transitive. In other words $A_x$  has the natural structure of a  complex affine space with associated (model) vector space $L(A)_x$.

The disjoint union $L(A)\edf \coprod_{x\in X} L(A)_x$
has the natural structure of a holomorphic vector bundle, which can be identified with the associated vector bundle $P_A\times_\rho \C^r$,  where $\rho:\Aff(r,\C)\to\GL(r,\C)$ is the standard linear representation   which assigns to each affine transformation its linear part. 
 
\begin{dt}
The holomorphic vector bundle $L(A)$ will be called the linearisation of $A$, and its associated locally free sheaf will be denoted by ${\cal L}(A)$.
\end{dt}

We will also need the representation $\tilde\rho:\Aff(r,\C)\to \GL(r+1,\C)$
of $\Aff(r,\C)$ defined by
$$
\tilde\rho(z\mapsto Mz+\zeta)\edf \begin{pmatrix}
1 &0\\ 
\zeta & M	
\end{pmatrix}.
$$
The image  $\tilde \Aff(r,\C)$  of this group monomorphism is a closed affine subgroup of $\GL(r+1,\C)$. The associated holomorphic vector bundle
$$
\tilde L(A)\edf P_A\times_{\tilde \rho}\C^{r+1}
$$
is a holomorphic rank $(r+1)$-vector bundle on $X$ which comes with a canonical $\tilde \Aff(r,\C)$-structure (reduction of structure group from $\GL(r+1,\C)$ to $\tilde \Aff(r,\C)$), see for instance \cite[Définition 2.3.1]{Te}. We will denote by $\tilde {\cal L}(A)$ the locally free sheaf associated with $\tilde L(A)$.

Note that 
\begin{re} \label{rem-p0-phi}
\begin{enumerate}
\item The subgroup 	$\tilde \Aff(r,\C)\subset \GL(r+1,\C)$ leaves invariant the projection $\C^{r+1}\ni (z_0,\dots,z_n)\stackrel{p_0}{\mapsto} z_0\in\C$.
\item Let $E$ be a holomorphic complex vector bundle of rank $r+1$ on $X$. The data of an 	$\tilde \Aff(r,\C)$-structure   on $E$ is equivalent to the data of a bundle epimorphism $\varphi:E\to X\times\C$. 
\end{enumerate}
\end{re}
Applying  this remark  to $\tilde L(A)$ we obtain a canonical bundle epimorphism
$$
\varphi_A:\tilde L(A)\to X\times\C 
$$
induced by the canonical 	$\tilde \Aff(r,\C)$-structure of $\tilde L(A)$. We will denote by the same symbol the associated sheaf epimorphism $\tilde {\cal L}(A)\to {\cal O}_X$.

\begin{dt}
An augmented rank $r+1$ holomorphic vector bundle on $X$ is a  	pair  $(E,\varphi)$ consisting of a holomorphic rank $r+1$-vector bundle $E$ on $X$  and a holomorphic bundle epimorphism $\varphi:E\to  X\times\C$ onto the trivial line bundle $X\times\C$ on $X$.
\end{dt}

Following \cite{Pl1}, \cite{Pl2} we will denote by $E_X^r$ the category of augmented rank $r+1$ holomorphic vector bundles on $X$.

\begin{re}\label{A-ag} Let $\ag=(E,\varphi)$ be an  augmented rank $r+1$ holomorphic vector bundle on $X$. Then the natural map $A_\ag\edf \varphi^{-1} (X\times\{1\})\to X$ comes with a natural structure ${\cal T}_\ag$ of a holomorphic affine bundle of rank $r$ on $X$ whose associated linearisation is canonically isomorphic to $\ker(\varphi)$.	
\end{re}

%
Moreover one has (see \cite[Lemma 2.2]{GoPr}, \cite{Pl1}, \cite{Pl2}):
\begin{pr} \label{AffBdls-AugmBdls}
The functors $A\mapsto (\tilde L(A),\varphi_A)$, $\ag\mapsto 	(A_\ag,{\cal T}_\ag)$ define an equivalence of groupoids between the groupoid  of holomorphic affine rank $r$-bundles on $X$ and the groupoid of augmented holomorphic rank $(r+1)$-vector bundles on $X$.

Via the first functor the linearisation $L(A)$ of an affine bundle $A$ is canonically isomorphic with $\ker(\varphi_A)$. 
\end{pr}

We will denote by $(e_A)$ the  obtained extension of locally free sheaves
$$
0\to \ker(\varphi_A)={\cal L}(A)\hookrightarrow \tilde{ {\cal L}}(A)\textmap{\varphi_A}  {\cal O}_X\to 0 \eqno{(e_A)}
$$
which will play an importent role in this article.

 Recall (see for instance \cite[Proposition 2 p. 184]{At}, \cite[Proposition 6.3 p. 234]{Ha}) that, in general, for any coherent sheaf ${\cal F}$ on $X$, isomorphism classes of extensions
 $$
 0\to {\cal F}\textmap{j} {\cal F}'\to {\cal O}_X\to 0 \eqno{(e)}
 $$
  of ${\cal O}_X$ by ${\cal F}$ correspond bijectively to elements of $\Ext^1({\cal O}_X,{\cal F})=H^1(X,{\cal F})$ via the extension class (invariant) map
  $$
  (e)\mapsto h_e\coloneq \delta(1),
  $$  
where $\delta:H^0(X,{\cal O}_X)\to H^1(X,{\cal F})$ is the connecting morphism associated with the exact sequence $(e)$.

\begin{dt}\label{DefhA}
	Let $A$ be an affine bundle over a complex manifold $X$ and $(e_A)$ the associated short exact sequence (see Section \ref{section-AffBundles}). We define the cohomology invariant of $A$ by  $h_A\edf h_{e_A}\in H^1(X,{\cal L}(A))$.
\end{dt}

We can now  state  the following set theoretical classification theorem (see \cite[Theorem 1.9]{Bou1}):

\begin{thr} \label{ClassTh}
Let $X$ be a complex manifold, $E$ a holomorphic vector bundle of rank $r$ on $X$, ${\cal E}$ the associated locally free sheaf and $\Aut(E)=\Aut({\cal E})$ its automorphism group. The map $A\mapsto h_A \in H^1(X,{\cal L}(A))$ induces a bijection  between the set of isomorphism classes of holomorphic affine  rank $r$-bundles $A$  on $X$ with $L(A)\simeq E$ and the quotient $H^1(X,{\cal E})/\Aut(E)$.
\end{thr}

%
%
%

By Cartan's finiteness theorem we have:

\begin{re} If $X$ is a  compact complex manifold, then the space $H^1(X,{\cal E})$ intervening in Theorem \ref{ClassTh} is a complex space of finite dimension. 
\end{re}

Moreover one has (see \cite[Proposition 2.5]{Bou1})

\begin{pr}\label{Aut(E)AffAlg}
Let $X$ be a compact connected complex manifold, $E$ a holomorphic vector bundle on $X$ and ${\cal E}$ the associated locally free sheaf. Then $\Aut(E)=\Aut({\cal E})$ is an affine algebraic group.
\end{pr}

\subsection{Gluing sheaves}\label{gluing-section}

We refer to \cite[sections II.1.2, II.1.6]{Go}  for the correspondence between sheaves of sets and étale spaces and the general formalism of sheaf homomorphisms.

We also refer to \cite[Exercise 1.22, p. 69]{Ha}, \cite[Section 1.6.33]{Stacks} for the following gluing principle for sheaves. In the following definition and proposition by sheaf we mean sheaf of sets or sheaf of abelian groups. 

Let $X$ be a topological space and  $\mathcal{U}=(U_i)_{i\in I}$  an open cover of $X$.

\begin{dt}
 A sheaf gluing data subordinate to 	$\mathcal{U}$ is a pair $\big(({\cal F}_i)_{i\in I},(\varphi_i^j)_{(i,j)\in I^2}\big)$ consisting of 
\begin{enumerate}
\item a family $({\cal F}_i)_{i\in I}$ of sheaves, ${\cal F}_i$ being a sheaf on $U_i$ for any $i\in I$.
\item a family $(\varphi_i^j)_{(i,j)\in I^2}$  of sheaf isomorphisms, $\varphi_i^j:{\cal F}_i|_{U_i\cap U_j}\textmap{\simeq}{\cal F}_j|_{U_i\cap U_j}$ satisfying the cocycle conditions:  
 \begin{enumerate}
	\item $\varphi_i^i=\id$  for any $i\in I$.
	\item For each $(i,j,k)\in I^3$  we have $\varphi_i^k=\varphi_j^k \circ \varphi_i^j$ on $U_i\cap U_j \cap U_k$.
\end{enumerate}
	\end{enumerate}

A morphism 
$$\big(({\cal F}_i)_{i\in I},(\varphi_i^j)_{(i,j)\in I^2}\big) \to \big(({\cal F}'_i)_{i\in I},({\varphi'}_i^j)_{(i,j)\in I^2}\big)$$
of gluing data subordinate to $\mathcal{U}$  is a family of morphisms $(u_i:{\cal F}_i\to {\cal F}'_i) _{i\in I}$   such that the diagram
$$
\begin{tikzcd}[row sep=6mm]
{\cal F}_i|_{U_i\cap U_j}\ar[d, "\varphi^j_i"']\ar[r, "u_i|_{U_i\cap U_j}"]&{\cal F}'_i|_{U_i\cap U_j}\ar[d, "\varphi'^j_i"]	\\
{\cal F}_j|_{U_i\cap U_j}\ar[r, "u_j|_{U_i\cap U_j}"]&{\cal F}'_j|_{U_i\cap U_j},\ 
\end{tikzcd}
$$
is commutative for any $(i,j)\in I^2$.
\end{dt}

We refer to \cite[Lemma 6.33.2]{Stacks} for the following result which will play a crucial role in our results: 
\begin{pr}\label{Gluing}
To any sheaf gluing data $\big(({\cal F}_i)_{i\in I},(\varphi_i^j)_{(i,j)\in I^2}\big)$ subordinate to ${\cal U}$ there is a canonically associated    pair $({\cal F},(\psi_i: \mathcal{F}_i \textmap{\simeq} \mathcal{F}|_{U_i} )_{i\in I})$, where ${\cal F}$ is a sheaf on $X$ and for any $i\in I$, $\psi_i: \mathcal{F}_i \textmap{\simeq} \mathcal{F}|_{U_i} $ is a sheaf isomorphism on $U_i$ which satisfies: 
\begin{enumerate}
\item The compatibility condition: 
$$\forall (i,j)\in I^2,\ \psi_i=\psi_j \circ \varphi_i^j \hbox{ on }U_i\cap U_j. \eqno{\rm (Co)}$$
\item The following universal property: 

For any pair $({\cal G},(g_i: \mathcal{F}_i \to  \mathcal{G}|_{U_i} )_{i\in I})$ consisting of a sheaf ${\cal G}$ on $X$ and a family of sheaf homomorphisms    satisfying the compatibility condition
$$\forall (i,j)\in I^2,\ g_i=g_j \circ \varphi_i^j \hbox{ on }U_i\cap U_j, \eqno{\rm (Co')}$$
there exists a  unique sheaf homomomorphism $g:{\cal F}\to  {\cal G}$ such that 
\begin{equation}\label{g_i-psi_i}
\forall i\in I,\ g|_{U_i}\circ \psi_i=g_i. 
\end{equation}
%
%
 
%
%
\end{enumerate} 
\end{pr}


%
\begin{dt}
The (essentially unique) sheaf given by Proposition 	\ref{Gluing} will be called to sheaf associated to the sheaf gluing data $\big(({\cal F}_i)_{i\in I},(\varphi_i^j)_{(i,j)\in I^2}\big)$ and will be denoted by $\coprod_{i\in I} {\cal F}_i \big/_{(\varphi_i^j)}$.
\end{dt}

\begin{re}\label{Mor-Ass-Mor-GluingData}
For any morphism  
$$(u_i:{\cal F}_i\to {\cal F}'_i) _{i\in I}:  \big(({\cal F}_i)_{i\in I},(\varphi_i^j)_{(i,j)\in I^2}\big) \to \big(({\cal F}'_i)_{i\in I},({\varphi'}_i^j)_{(i,j)\in I^2}\big)$$
	of gluing data, there exists a unique morphism 
	$$\coprod_{i\in I} {\cal F}_i \big/_{(\varphi_i^j)}\eqcolon  {\cal F}\textmap{\coprod_{i\in I} u_i} {\cal F}'\coloneq \coprod_{i\in I} {\cal F}'_i \big/_{(\varphi_i^j)}$$ 
such that the following diagram is commutative
$$
\begin{tikzcd}[row sep=5mm, column sep=1.8cm]
{\cal F}_i\ar[d, "\psi_i"']\ar[r, "u_i"]&{\cal F}'_i\ar[d, "{\psi'}_i"]	\\
{\cal F}\big|_{U_i}\ar[r, "\coprod_{i\in I} u_i|_{U_i}"]& {\cal F}'\big|_{U_i}\ 
\end{tikzcd}, 
$$
where  $({\cal F},(\psi_i: \mathcal{F}_i \textmap{\simeq} \mathcal{F}|_{U_i} )_{i\in I})$, $({\cal F}',(\psi'_i: \mathcal{F}'_i \textmap{\simeq} \mathcal{F}'|_{U_i} )_{i\in I})$ are the pairs associated with the gluing data $\big(({\cal F}_i)_{i\in I},(\varphi_i^j)_{(i,j)\in I^2}\big)$, $ \big(({\cal F}'_i)_{i\in I},({\varphi'}_i^j)_{(i,j)\in I^2}\big)$.
\end{re}
\begin{re}\label{Glue-sheaves-Ox-modules}
Let $(X,{\cal R})$ be a ringed space (complex space). Proposition \ref{Gluing} remains true in the category of (respectively coherent)  ${\cal R}$-modules.

\end{re} 

We will need the following obvious properties of the sheaves and morphisms associated with gluing data (see \cite[Lemma 6.33.4]{Stacks} ) 
\begin{re}\label{FunctorialityGluing}
Let $(X,{\cal R})$ be a ringed space and ${\cal U}=(U_i)_{i\in I}$ an open cover of $X$. The assignment 	
$$
\big(({\cal F}_i)_{i\in I},(\varphi_i^j)_{(i,j)\in I^2}\big)\mapsto \coprod_{i\in I} {\cal F}_i \big/_{(\varphi_i^j)}
$$
has the following properties:
\begin{enumerate}
\item 	Is an exact functor on the category of ${\cal R}$-module gluing data in the following sense: let
$$
\big(({\cal F}_i)_{i\in I},(\varphi_i^j)_{(i,j)\in I^2}\big),\ \big(({\cal F}'_i)_{i\in I},(\varphi'^j_i)_{(i,j)\in I^2}\big),\ \big(({\cal F}_i)_{i\in I},({\varphi''}_i^j)_{(i,j)\in I^2}\big)
$$
be ${\cal R}$-module gluing data and let $({\cal F}'_i\textmap{u_i} {\cal F}_i)_{i\in I}$ $({\cal F}_i\textmap{v_i} {\cal F}''_i)_{i\in I}$  be morphisms of gluing data such that ${\cal F}'_i\textmap{u_i}{\cal F}'\textmap{v_i} {\cal F}''_i$ is exact for any $i\in I$.
    Then the sequence
$$
\coprod_{i\in I} {\cal F}'_i \big/_{({\varphi'}_i^j)}\textmap{\coprod_{i\in I} u_i}\coprod_{i\in I} {\cal F}_i \big/_{(\varphi_i^j)}\textmap{\coprod_{i\in I} v_i}\coprod_{i\in I} {\cal F}''_i \big/_{({\varphi''}_i^j)}
$$
is exact.
\item Commutes with tensor products in the following sense:

Let $\big(({\cal F}_i)_{i\in I},(\varphi_i^j)_{(i,j)\in I^2}\big)$ be an ${\cal R}$-module gluing data, 
$$({\cal F}=\coprod_{i\in I} {\cal F}_i \big/_{(\varphi_i^j)},(\psi_i: \mathcal{F}_i \textmap{\simeq} \mathcal{F}|_{U_i} )_{i\in I})$$
 the associated pair  and ${\cal A}$  an ${\cal R}$-module on $X$. Then the family of ${\cal R}$-module isomorphisms
$\big(\psi_i\otimes\id_{\cal A}:{\cal F}_i\otimes{\cal A}|_{U_i}\to (\mathcal{F}\otimes {\cal A})|_{U_i}\big)_{i\in I}$
induces a canonical isomorphism 
$\coprod_{i\in I} ({\cal F}_i\otimes{\cal A}|_{U_i} )\big/_{(\varphi_i^j\otimes\id_{\cal A})}\textmap{\simeq}\big(\coprod_{i\in I} {\cal F}_i \big/_{(\varphi_i^j)}\big)\otimes{\cal A}$.
 \end{enumerate}

\end{re}

\subsection{Enlarged extensions}\label{EnlargedExtSection}

Let $X$ be a topological space, $\iota:{\cal A}\to {\cal A}'$ a monomorphism of sheaves of abelian groups on $X$, and  
$$
0\to {\cal A}\textmap{u}{\cal B}\textmap{v} {\cal C}\to 0 \eqno{(e)}
$$
an exact sequence  of sheaves of abelian groups on $X$, regarded as an extension of ${\cal C}$ by ${\cal A}$. 
\begin{dt} \label{iota-enlarged-def}
	The enlarged extension associated with $(e)$ and $\iota$  is 
	$$
	0\to {\cal A}'\textmap{u'}{\cal B}'\textmap{v'} {\cal C}\to 0 \eqno{(e')}
	$$
	where 
	\begin{enumerate}
	\item $
	{\cal B}'\edf \coker \big( (\iota, -u) : {\cal A}\to {\cal A}'\oplus {\cal B}\big)$,
	\item $u'$ is the   composition 
	$$
{\cal A}'\textmap{i_1}{\cal A}'\oplus {\cal B}\to \qmod{{\cal A}'\oplus {\cal B}}{\im\big((\iota, -u)\big)}\eqcolon{\cal B}',$$
\item $v'$ is induced by the composition ${\cal A}'\oplus {\cal B}\textmap{p_2}{\cal B}\textmap{v}{\cal C}$
whose kernel obviously contains $\im(\iota, -u)$.
	
	\end{enumerate}
 
\end{dt}

Note that the kernel of the composition 
${\cal A}'\oplus {\cal B}\textmap{p_1}{\cal A}\map \qmod{{\cal A'}}{\iota({\cal A})}=\coker(\iota)$
contains $\im(\iota, -u)$, so this composition induces an epimorphism $q_e^\iota:{\cal B}'\to \coker(\iota)$ whose kernel is $\im(b')$ and will be called {\it the canonical epimorphism} on the central term ${\cal B}'$ of the enlarged extension. 
\begin{re}\label{enlarged-remark1}
The enlarged extension $(e')$ defined above has the following properties:
\begin{enumerate}
\item Fits in  the commutative diagram 
$$
\begin{tikzcd}[row sep=5mm]
&0\ar[d]&0\ar[d]&0\ar[d]&\\
0\ar[r]&{\cal A}\ar[r,"u"]\ar[d, "\iota"]&{\cal B}\ar[r,"v"]\ar[d,"b' "]&{\cal C}\ar[r]\ar[d, equal] &0\\
0\ar[r] &{\cal A}'\ar[r,"u'"]\ar[d] &{\cal B}'\ar[r,"v' "]\ar[d, "q_e^\iota"]&{\cal C}\ar[r] &0\\
&\coker(\iota)\ar[r,equal]\ar[d]&\coker(\iota)\ar[d]&&\\
&0&0&&	
\end{tikzcd}
$$
with exact rows and columns.	 The vertical arrow $q_e^\iota$ is the unique morphism ${\cal B}'\to \coker(\iota)$ which makes the diagram commutative and the second column exact.
\item Satisfies the following universal property: for any commutative diagram
$$
\begin{tikzcd}
0\ar[r]&{\cal A}\ar[r,"u"]\ar[d, "\iota"]&{\cal B}\ar[r,"v"]\ar[d,"b''"]&{\cal C}\ar[r]\ar[d, equal] &0\\
0\ar[r] &{\cal A}'\ar[r,"u''"]  &{\cal B}''\ar[r,"v'' "] &{\cal C}\ar[r] &0	
\end{tikzcd}
$$
with exact rows there a unique isomorphism $\beta:{\cal B}'\textmap{\simeq} {\cal B}''$ making commutative the following diagram
$$
\begin{tikzcd}[column sep=1.7cm, row sep=6mm]
0\ar[r]&{\cal A}\ar[r,"u"]\ar[d, "\iota"]&{\cal B}\ar[r,"v"]\ar[d,  "b" description]\ar[dd,  bend right =38, near start, "b''"' description]&{\cal C}\ar[r]\ar[d, equal] &0\\
0\ar[r]&{\cal A}'\ar[r,crossing over, near start, "u'"]\ar[d, equal]&{\cal B}'\ar[r,"v'"]\ar[d, dashed, "\beta" description, "\simeq"]&{\cal C}\ar[r]\ar[d, equal] &0\\
0\ar[r] &{\cal A}'\ar[r,"u''"]  &{\cal B}''\ar[r,"v'' "] &{\cal C}\ar[r] &0	
\end{tikzcd}
$$
	
\end{enumerate}

\end{re}

\begin{proof} 
(1) Direct checking.

(2) The unique  isomorphism $\beta:{\cal B}'\to {\cal B}''$  which makes the above diagram commutative is the isomorphism induced by $(u'',b''):{\cal A}'\oplus {\cal B}\to {\cal B}''$. \end{proof}

\begin{re}\label{enlarged-remark2} Let $(X,{\cal R})$ be a ringed space (complex space).
The definition of the enlarged extension and its properties stated in Remark \ref{enlarged-remark1} remain true in the category of sheaves of (respectively coherent) ${\cal R}$-modules.  
\end{re}

In section \ref{Equiv-Of-Stacks} we will need the following special class of enlarged extensions: 

Let $X$ be a complex manifold and $\varphi:\tilde {\cal E}\to {\cal O}_X$  be a surjective morphism of locally free sheaves of $X$ such that $\mathcal{E}_\varphi\edf\ker(\varphi)\simeq\mathcal{E}$, where $\mathcal{E}$ is the locally free sheaf associated with a vector bundle $E$ on $X$. We have an exact sequence

$$
0\to \mathcal{E}_\varphi \xhookrightarrow{} \tilde{\mathcal{E}} \textmap{\varphi} \mathcal{O}_X \to 0, \eqno{(e_\varphi)}
$$
which will be regarded as an extension of ${\cal O}_X$ by ${\cal E}_\varphi$. 
Since ${\cal E}_\varphi$ is locally free, we also have the exact sequence
$$
0\to \mathcal{E}_\varphi \xhookrightarrow{\iota_D^\varphi} \mathcal{E}_\varphi(D)\textmap{\pi_D^\varphi} \mathcal{E}_\varphi(D)_D\to 0 \eqno{(c_\varphi^D)}.
$$
obtained from the canonical exact sequence $0\to \mathcal{O}_X \xhookrightarrow{} \mathcal{O}_X(D)\to \mathcal{O}_D(D)\to 0$ by tensorising with ${\cal E}_\varphi$.

\begin{dt}\label{enlarged-ext} We define  the $D$-enlarged extension associated with $(e_\varphi)$ to be the enlarged extension
$$
0 \to \mathcal{E}_\varphi(D) \textmap{\iota^D_\varphi} \tilde{\mathcal{E}}_\varphi^D \textmap{\varphi^D} \mathcal{O}_X\to 0. \eqno{(e_\varphi^D)}
$$
 associated to $(e_\varphi)$ and the canonical sheaf monomorphism ${\cal E}_\varphi\hookrightarrow {\cal E}_\varphi(D)$ in the sense of Definition \ref{iota-enlarged-def}.

In other words  $\tilde{\mathcal{E}}_\varphi^D$ is the cokernel of the morphism $(\iota_D^\varphi ,-j):\mathcal{E}_\varphi \to \mathcal{E}_\varphi(D)\oplus \tilde{\mathcal{E}}$, and $\tilde{\mathcal{E}}_\varphi^D \textmap{\varphi^D} \mathcal{O}_X$ is the epimorphism induced by $\varphi$ whose kernel is obviously identified with $\mathcal{E}_\varphi(D)$. 

The canonical  monomorphism into  the central term of $\tilde{\mathcal{E}}_\varphi^D$ of $(e_\varphi^D)$ is the the sheaf monomorphism given by the composition 
$$
\begin{tikzcd}
\tilde {\cal E} \ar[r, hook] \ar[rr, bend right =20, "j^D_\varphi"']& {\cal E}_{\varphi}(D)\oplus \tilde {\cal E} \ar[r, two heads] &	\tilde{\mathcal{E}}_\varphi^D.
\end{tikzcd}
$$

The canonical epimorphism on   $\tilde{\mathcal{E}}_\varphi^D$ of $(e_\varphi^D)$ is the epimorphism 
$$
q^\varphi_D:\tilde{\mathcal{E}}_\varphi^D\to {\cal E}_\varphi(D)_D
$$
induced by the composition 
$$\mathcal{E}_\varphi(D)\oplus \tilde{\mathcal{E}}\textmap{pr_1}   \mathcal{E}_\varphi(D)\textmap{\pi_D^\varphi} \mathcal{E}_\varphi(D)_D.$$
\end{dt}

Taking into account the general properties of the enlarged extension (see section \ref{EnlargedExtSection}), we have:
\begin{re} \label{PropOfCanonEpim}
The extensions $(e_\varphi)$,  $(e_\varphi^D)$ fit in the commutative diagram
\begin{equation}\label{restriction-diag}
\begin{tikzcd}[row sep=5mm]
&0\ar[d] &0 \ar[d] \\[-6pt]
0 \ar[r] & {\cal E}_\varphi  \ar[r, hook]  \ar[d, hook'] & \tilde {\cal E}  \ar[r, "\varphi" ]\ar[d, "j^D_\varphi"'] & {\cal O}_X \ar[r]\ar[d, equal]& 0 \\
0  \ar[r]  & \mathcal{E}_\varphi(D)  \ar[d] \ar[r, "\iota^D_\varphi"] &  \tilde{\mathcal{E}}_\varphi^D \ar[d, "q_D^\varphi"] \ar[r, "\varphi^D"'] & \mathcal{O}_X \ar[l,"\sigma"', dashed ,bend right =40] \ar[r] & 0 \\
 & {\cal E}_\varphi(D)_D \ar[r, equal]\ar[d]& {\cal E}_\varphi(D)_D \ar[d] & & \\[-6pt]
 & 0 & 0
\end{tikzcd}
\end{equation}
with exact rows and columns. The canonical epimorphism $q_D^\varphi : \tilde{\mathcal{E}}_\varphi^D \to {\cal E}_\varphi(D)_D$ is  the unique morphism 	$\tilde{\mathcal{E}}_\varphi^D \to {\cal E}_\varphi(D)_D$ which makes the second column exact and the diagram commutative.

\end{re}

Note that 

\begin{re}\label{Extension-Invariant-Remark}
The extension invariant of $(e_\varphi^D)$ is precisely the image of the extension invariant $h_\varphi\in H^1(X,{\cal E}_\varphi)$ of $(e_\varphi)$ in  $H^1(X,{\cal E}_\varphi(D))$ via the morphism $H^1(\iota^\varphi_D)$.	
\end{re}

\subsection{The representability of the functor  \texorpdfstring{$S\mapsto \Aut_{S\times X}(S\times E)$}{aut(E)}}  {\ }

Let  $f:X\to S$ be a morphism of complex spaces and ${\cal F}$ be a coherent sheaf on  $X$. For $s\in S$ we will denote by $X_s$ the fibre of $f$ over $s$ and by ${\cal F}_s$ the restriction of ${\cal F}$ to $X_s$. Now let
$$
\begin{tikzcd}[row sep=5mm]
     X'   \ar[d, "f’"]\ar[r, "G"]  & X \ar[d, "f"] \\
     S' \ar[r, "g"] & S 
\end{tikzcd}
$$
be a commutative diagram of complex spaces, and $\mathcal{F}$ be a coherent $\mathcal{O}_X$-module.

Recall (see for instance in \cite[p. 93]{BS}) that for any $q\in \N$ we have a canonical  morphism of sheaves: 
$$
g^*(R^qf_{*}({\cal F}))\to R^q (f')_{*}(G^*({\cal F})).
$$
If $f:X \to S$ is proper, $X'\edf S'\times_S X$, ${\cal F}$ is flat over $S$  and cohomologically flat in dimension $q$, this canonical morphism is in fact an isomorphism (see \cite[Theorem 4.10 (d), p. 133]{BS}). 	

On the other hand, by Grauert's Theorem (see \cite[Theorem 4.12 (ii), p. 134]{BS}), we know that if moreover 
\begin{enumerate}
	\item\label{remC1} $S$ is reduced, and, 
	\item\label{remC4} the map $S \ni s \mapsto h^q(X_s,{\cal F}_s)$ is locally constant,
\end{enumerate}
then   ${\cal F}$ is cohomologically flat in dimension $q$. We now apply this result to the diagram 
$$
\begin{tikzcd}[row sep=5mm]
    S \times X   \ar[d, "p_{S}"'] \ar[r, "p_X"]  & X \ar[d] \\
     S \ar[r] & \{\star\}
\end{tikzcd}
$$
where $\{\star\}$ is a reduced singleton and $q=0$. Since the conditions (\ref{remC1}) and (\ref{remC4}) above are obviously satisfied, we obtain the following:  

\begin{re}\label{pSpX}

Let $X$ be a compact complex space and $\mathcal{F}$ a coherent sheaf on $X$. For any complex space $S$, the canonical morphism 
$$
\mathcal{O}_S \otimes_{\C} H^0(X,{\cal F})\to (p_S)_{*}(p_X^*(\mathcal{F}))
$$
is an isomorphism.  Applying the functor $H^0(S,-)$ to this sheaf isomorphism, we obtain a canonical isomorphism 
\begin{equation}\label{can-isos}
	{\cal O}(S,H^0(X,\mathcal{F})) \to H^0(S,(p_S)_{*}((p_X^*)(\mathcal{F})))=H^0(S\times X, (p_X^*)(\mathcal{F})).
\end{equation}
\end{re}

Let now $X$ be a compact complex manifold, $E$ a holomorphic vector bundle on $X$ and ${\cal E}$ its associated locally free sheaf. Recall (see Proposition \ref{Aut(E)AffAlg}) that the automorphism group $\Aut(E)=\Aut({\cal E})$ is an affine algebraic group. From now on we will denote by $G$ the group $\Aut(E)$ endowed with its natural affine algebraic group structure given by this proposition.

We have the following  
\begin{pr}\label{repres-Aut}
	The functor $S\mapsto \Aut_{S\times X}((p_X^{S\times X})^*(E))$ is representable by $G$.
\end{pr}
\begin{proof}

Using the isomorphism (\ref{can-isos}) when $\mathcal{F}=\mathcal{E} nd(\mathcal{E})=\mathcal{E}^{\vee}\otimes \mathcal{E}$, we obtain identifications
\begin{equation}\label{iso-End-H0}
\begin{split}
	\mathcal{O}_S(\End(\mathcal{E}))(S) \textmap{\simeq } H^0(S\times X, p_X ^*(\mathcal{E}^\vee \otimes  \mathcal{E}))&=H^0(S\times X, p_X ^*(\mathcal{E})^\vee \otimes p_X^*(\mathcal{E}))
	\\&=\End(p_X^*(\mathcal{E})),
\end{split}
\end{equation}
which are functorial with respect to $S$ and compatible with compositions. Here $p_X=p_X^{S \times X}$ stands for the projection on the first factor. Now note that, since $G=\Aut(E)$ is open in $\End(E)$, the set ${\cal O}_S(G)(S)$ of complex space morphisms $S\to G$ coincides with the  subset of ${\cal O}_S(\End(E))(S)$ consisting of complex space morphisms $S\to \End(E)$ whose underlying map takes values in  $G$. On the other hand the set of invertible elements in ${\cal O}_S(\End(E))(S)$ is obviously ${\cal O}_S(G)(S)$. Taking the set of invertible elements in (\ref{iso-End-H0}), we obtain a group isomorphism  
$$
{\cal O}_S(G)(S)\textmap{\simeq}\Aut_{S\times X}(S\times E)=\Aut_{S\times X}(p_X^*(\mathcal{E}))
$$
which is functorial  with respect to $S$.
\end{proof}
\begin{re} The analogue of the representability theorem given by Proposition \ref{repres-Aut} in the algebraic geometric category is known for  principal $G$-bundles over proper algebraic varieties, where $G$ is an arbitrary algebraic group (see \cite[Theorem 4.2 p. 16, comments on p. 17]{Br}). \end{re}
\subsection{Fibrewise constant families of vector bundles}\label{fbrwse-isom-E-sect}

Let $X$ be a compact complex manifold and $S$ a complex space. For a holomorphic vector bundle $F$ (locally free sheaf ${\cal F}$) on $S\times X$ we will denote by $F_s$ (respectively $\mathcal{F}_s$) the restriction  $F|_{\{s\}\times X}$ (respectively $\mathcal{F}|_{\{s\}\times X}$) regarded as a vector bundle (respectively locally free sheaf) on $X$.
Let  $E$ be a  holomorphic vector bundle of rank $r$ on $X$ and $\mathcal{E}$ its associated locally free sheaf.
%
%
%
\begin{dt}\label{fibrewise-dt}
	 A vector bundle $F$ over $S\times X$ is said to be fibrewise isomorphic to $E$ if for any $s\in S$ there exists an open neighbourhood  $U$ and an isomorphism
	$$
	 (p_X^{U\times X})^*(E)\textmap{\simeq} F|_{U\times X} ,
	$$
of holomorphic vector bundles 	over $U\times X$, where $p_X^{U\times X}:U\times X\to X$ stands for the projection on the second factor.

Similarly, a locally free sheaf ${\cal F}$ over $S\times X$ is said to be fibrewise isomorphic to ${\cal E}$ if for any $s\in S$ there exists an open neighbourhood  $U$ and an isomorphism $(p_X^{U\times X})^*({\cal E})\textmap{\simeq} {\cal F}|_{U\times X}$
	of ${\cal O}_{U\times X}$-modules.
\end{dt}

%

%

\begin{re}
	Let $F \to S\times X$ be a vector bundle fibrewise isomorphic to $E$. For each $s\in S$, the vector bundle $F_s$ is isomorphic to $E$. 
\end{re}

The converse of this implication holds when $S$ is reduced, so for reduced space of parameters $S$ the condition ``$F$ is fibrewise isomorphic to $E$" agrees with the intuitive interpretation of this condition. More precisely
\begin{pr}
	Let $S$ be a \textbf{\textit{reduced}} complex space and $F$ an holomorphic bundle on $S \times X$ such that for any $s\in S$ the vector bundle $F_s$ is isomorphic to $E$. Then $F$ is fibrewise isomorphic to $E$ in the sense of Definition \ref{fibrewise-dt}.  
\end{pr}
\begin{proof}
Let $\mathcal{H}$ be the locally free sheaf $p_X^*(\mathcal{E})^\vee\otimes \mathcal{F}$ over $S \times X$. Since $p_S: S \times X \longrightarrow S$ is flat, $\mathcal{H} $ is flat over $S$. Since $X$ is compact, the projection $p_S$ is   proper. We have:
$$
\mathcal{H}_s \simeq  \mathcal{E}^\vee\otimes \mathcal{F}_s  \simeq \mathcal{E}^\vee\otimes \mathcal{E} .
$$
Therefore the map $S\ni s\mapsto h^0(X_s, \mathcal{H}_s)$ is constant and equals  $d\edf \dim_\C(\End(\mathcal{E}))$. Since $S$ is reduced, it follows  by Grauert's Theorem \cite[Theorem 4.12 (ii)]{BS}, that $\mathcal{H}$ is cohomologically flat in dimension 0 over $S$ \cite[p. 133]{BS}, and in particular:

\begin{enumerate}
\item  the sheaf $p_{S_{*}}(\mathcal{H})$ is locally free of rank $d$, and 	
\item for each $s\in S$ the natural linear map 
$$
\qmod{(p_S)_{*}(\mathcal{H})_s}{\mg_s (p_S)_{*}(\mathcal{H})_s}\edf (p_S)_{*}(\mathcal{H})(s)\to H^0(X_s ,\mathcal{E}^\vee\otimes \mathcal{F}_s)\simeq \Hom(\mathcal{E},\mathcal{F}_s)
$$
is an isomorphism (see \cite[Theorem 4.10 (c)]{BS}).
\end{enumerate}
Let $s\in S$.  Choose an isomorphism $\alpha\in  \Hom(\mathcal{E},\mathcal{F}_s)$. Such an isomorphism always exists, because by assumption $\mathcal{F}_s\simeq \mathcal{E}$. Let $\tilde{\alpha}\in (p_S)_{*}(\mathcal{H})_s$ be an element which is mapped to $\alpha$ via the composition 
$$
(p_S)_{*}(\mathcal{H})_s\to \qmod{(p_S)_{*}(\mathcal{H})_s}{\mg_s (p_S)_{*}(\mathcal{H})_s}\simeq \Hom(\mathcal{E},\mathcal{F}_s),
$$
and let $U$ be an open neighbourhood of $s$ such that $\tilde{\alpha}$ extends to an element of $(p_S)_{*}(\mathcal{H})(U)=H^0(U\times X,\mathcal{H})$.  Taking $U$ sufficiently small, we obtain an isomorphism $p_X^*({\cal E})|_{U\times X}\textmap{\simeq} {\cal F}|_{U\times X}$.
\end{proof}

Using similar arguments, one obtains the following result which gives an equivalent definition for the condition ``fibrewise isomorphic to $E$":
\begin{pr}
Let $S$ be a complex space and $F$ a vector bundle on $S\times X$. The following conditions are equivalent:
\begin{enumerate}
	\item $F$ is fibrewise isomorphic to $E$.
	\item For any $s\in S$, $F_s$ is isomorphic to $E$ and the locally free sheaf $p_X^*({\cal E})^\vee \otimes {\cal F}$ on $S\times X$ is cohomogically flat in dimension 0.
\end{enumerate} 	
\end{pr}

\begin{pr}\label{the-torsor-T^F}
Let $F$ be a holomorphic bundle on $S\times X$ fibrewise isomorphic to $E$. The sheaf of sets ${\cal T}^F$ on $S$ defined by 
\begin{equation*}
\begin{split}
{\cal T}^F(U)\edf \big\{ (p_X^{U\times X})^*(E)\textmap{h\simeq}F|_{U\times X}|\ h \hbox{ isomorphism of    }& \hbox{holomorphic vector} \\ &\hbox{ bundles over }U\times X\big\}
\end{split}
\end{equation*}	
is naturally an ${\cal O}_S(G)$-torsor. 
\end{pr}
\begin{proof} 
Note first for every open set $U$ of $S$, we have a natural right $\mathcal{O}_S(G)(U)$-action 
$$
{\cal T}^F(U)\times \mathcal{O}_S(G)(U)\to  {\cal T}^F(U).
$$
This action is free and transitive. Indeed, for any pair $(h,k)\in {\cal T}^F(U)\times{\cal T}^F(U)$ there exist a unique $g\in \mathcal{O}_S(G)(U)$ such that $k=h\circ g$, namely the element corresponding to $h^{-1}\circ k\in \Aut_{U\times X}((p_X^{U\times X})^*(E))$ under the isomorphism given by the representability Proposition \ref{repres-Aut}.

 On the other hand for any $s\in S$ the stalk ${\cal T}^F_s$ is non-empty because $F$ is fibrewise isomorphic to $E$. Therefore ${\cal T}^F_s$ is a ${\cal O}_S(G)_s$-torsor, which completes the proof.   
\end{proof}

Conversely, for any $\mathcal{O}_S(G)$-torsor ${\cal T}$, we can associate functorially an holomorphic vector bundle $F^{\cal T}\to S\times X$ which is fibrewise isomorphic to $E$ which is given explicitly by:  
$$
F^{\cal T}\edf \bigg(\coprod_{\tau \scriptsize  \hbox{    local section of }{\cal T} }U_\tau \times E \bigg)\bigg/ \sim, 
$$
where $\sim$ stands for the equivalence relation generated by the holomorphic vector bundle isomorphisms 
$$
\begin{tikzcd}
 (U_\tau\cap U_\theta)\times E \ar[d, hook']&  \ar[l, "g_{\tau\theta}"' ] (U_\tau\cap U_\theta)\times E\ar[d, hook'] \\
 U_\tau\times E & U_\theta\times E
 \end{tikzcd}
$$ 
where $g_{\tau \theta}\in {\cal O}_S(G)(U_{\tau} \cap U_{\theta})$ is defined by the condition 
$$
\theta|_{U_\tau\cap U_\theta}=\tau|_{U_\tau\cap U_\theta}\ g_{\tau\theta}. 
$$

\begin{pr}\label{composition-equiv}
	Let $F\to S \times X $ be a holomorphic vector bundle fibrewise isomorphic to $E$ and $\mathcal{T}$ a $\mathcal{O}_S(G)$-torsor on $S$. We have a natural isomorphisms 
	$$F\simeq F^{{\cal T}^F}, \ \mathcal{T} \simeq \mathcal{T}^{F^\mathcal{T}}.$$
	
	In particular, the functors $F \to \mathcal{T}^F$ and $\mathcal{T} \to F^\mathcal{T}$ define an equivalence between the groupoid of holomorphic vector bundles on $S \times X$ which are fibrewise isomorphic to $E$ and the groupoid of  $\mathcal{O}_S(G)$-torsors on $S$. 
	
	In particular, taking into account the well-known equivalence between principal $G$-bundles and $\mathcal{O}_S(G)$-torsors, we obtain  an equivalence between the groupoid of holomorphic vector bundles on $S \times X$ which are fibrewise isomorphic to $E$ and the groupoid of  principal $G$-bundles on $S$. 	
\end{pr}	                                                                                                                         
	    Using the formalism of associated  bundles, we have                                                                                                                      
\begin{re}
	Let $P \to S$ be a principal $G$-bundle over $S$. The holomorphic vector bundle on $S\times X$ fibrewise isomorphic to $E$ associated to $P$ (via the equivalence given by Proposition \ref{composition-equiv}) is precisely the associated bundle $P\times_\rho E$, regarded as a bundle over $S \times X$, where $\rho : G \times E \to E$ is the natural left action. 
\end{re}
\section{An equivalence of stacks}\label{Equiv-Of-Stacks}

%

Let $X$ be a compact complex manifold, $E$ be a holomorphic vector bundle of rank $r$ on $X$, ${\cal E}$ the associated locally free sheaf. We will assume that the pair $(X,{\cal E})$ satisfies:
\begin{main-assumpt} 
$X$ admits an effective divisor $D$ such that $$H^1(X,{\cal E}(D))=0.$$
\end{main-assumpt}

Note that this condition is automatically satisfied when $X$ is a projective complex manifold. In this case one can consider an ample effective divisor $H\subset X$ and, making use of Cartan-Serre-Grothedieck's vanishing theorem,  set $D=nH$ for sufficiently large $n$.  

\vspace{2mm}

From now on
\begin{itemize}
\item 	We  suppose that $(X,{\cal E})$ satisfies the main assumption,
\item We  fix  an effective divisor $D$ on $X$ such that $H^1(X,{\cal E}(D))=0$.
\item  We put $G\edf \Aut({\cal E})$.
\end{itemize}

For $g\in G$ we will denote by $g(D)$,  $g(D)_D$ the automorphisms induced by $g$ on the sheaves $\mathcal{E}(D)$, respectively $\mathcal{E}(D)_D$ and also on their spaces of global sections. 

\begin{dt} We define
 $$G\edf \Aut(\mathcal{E}), \ \ V_D=H^0(X,\mathcal{E}(D)_D), \ \ W_D\edf H^0(X,{\cal E}(D)),\ G_D\edf W_D\rtimes G,$$ 
where the semidirect product giving  $G_D$ is associated with the group morphism
$$
 G \ni g  \longmapsto H^0(g(D))\in \GL(H^0(X,{\cal E}(D)))=\GL(W_D).
$$	
\end{dt}
Taking into account the definition of semidirect product, the group law of $G_D$ is given by 
$$
(v,g)  (w,h)=(v+g(D)(w),g\circ h). 
$$

We also define 	a natural left $G_D$-action $\alpha_D:G_D\times V_D\to V_D$ on $V_D$ by
$$
\alpha_D((w,g), \beta) =g(D)_{D}(\beta)- w_D.
$$
where $w_D$ is the image of $w$ under the obvious  restriction morphism 
$$H^0(X,\mathcal{E}(D))=W_D\to V_D=H^0(X,\mathcal{E}(D)_D).$$

In this section we will  show that the stack $\mathscr{A}{\rm ff}_E$ defined in the introduction is isomorphic to the quotient stack $[V_D/G_D]$.

Let $S$ be a complex space.  As explained in the introduction,  we denote by $\tensor[_S]{\Afff}{_E}$  the fibre of $\mathscr{A}{\rm ff}_E$ over $S$, i.e. the groupoid whose objects are  augmented bundles $\tilde {\cal E}\textmap{\varphi} {\cal O}_{S\times X}$ on $S\times X$ such that ${\cal E}_\varphi\edf \ker(\varphi)$ is fibrewise isomorphic to ${\cal E}$ in the sense of Definition \ref{fibrewise-dt}. For such  augmented bundle $\tilde {\cal E}\textmap{\varphi} {\cal O}_{S\times X}$ over $S\times X$ and for an open subspace $U\subset S$, we will denote by $\tilde{\cal E}_U\textmap{\varphi_U}{\cal O}_{U\times X}$ its restriction to $U\times X$.

By Proposition \ref{the-torsor-T^F}, the sheaf of sets ${\cal T}_\varphi$ on $S$ defined by 
$$
{\cal T}_\varphi(U)\edf \big\{ (p_X^{U\times X})^*({\cal E})\textmap{f\simeq}{\cal E}_{\varphi_U}|\ f \hbox{ is an isomorphism of ${\cal O}_{U\times X}$-modules}\big\}
$$	
(for open sets $U\subset S$) is naturally an ${\cal O}_S(G)$-torsor. 
 
We will also denote by  $\tensor[_S]{\Torss}{}^{\alpha_D}_{G_D}$ the groupoid whose objects are pairs $(\Theta,\rho)$ consisting of an ${\cal O}_S(G_D)$-torsor  on $S$ and a ${\cal O}_S(G_D)$-equivariant sheaf morphism $\rho:\Theta\to {\cal O}_S(V_D)$. 

Note that  $\tensor[_S]{\Torss}{}^{\alpha_D}_{G_D} $ is precisely the fibre over $S$ of the quotient stack $[V_D/G_D]$ regarded as category fibred in groupoids. 

We will define an equivalence of categories using the functors 
$$
\tensor[_S]{\Tg}{}: \tensor[_S]{\Afff}{}_{E}\to \tensor[_S]{\Torss}{}^{\alpha_D}_{G_D},\ \tensor[_S]{\Ag}{}: \tensor[_S]{\Torss}{}^{\alpha_D}_{G_D}\to \tensor[_S]{\Afff}{}_{E}
$$
which will be defined in the following sections (sections \ref{functorTGS}, \ref{functorAGS}).

\subsection{The functor \texorpdfstring{$\tensor[_S]{\Tg}{}: \tensor[_S]{\Afff}{}_{E}\to \tensor[_S]{\Torss}{}^{\alpha_D}_{G_D}$}{TgS}} 
\label{functorTGS}
{\ }\\

Put $D_S\edf S\times D$ regarded as a divisor on $S\times X$. 
\begin{lm}\label{R1=0}
We have $R^1p_{S*}({\cal E}_\varphi(D_S))=0$.	
\end{lm}
\begin{proof}
Indeed, since ${\cal E}_\varphi$ is locally (with respect to $S$) isomorphic to $p_X^*({\cal E})$, it follows that ${\cal E}_\varphi(D_S)$ is locally (with respect to $S$) isomorphic to $p_X^*({\cal E}(D))$. The claim follows using the condition $H^1(X,{\cal E}(D))=0$ (our Main Assumption) and  the base change theorem applied to the morphism $X\to \{*\}$.
\end{proof}

Let $\tilde{\cal E}\textmap{\varphi}{\cal O}_{S\times X}$ be an object of the groupoid $\Afffs_E$. Consider the extension
$$
0\to {\cal E}_\varphi\hookrightarrow \tilde {\cal E}\textmap{\varphi}{\cal O}_{S\times X}\to 0 
$$
and its $D_S$-enlarged extension  
\begin{equation}\label{enl-ext-on-S}
0\to {\cal E}_\varphi(D_S)\textmap{\iota^D_\varphi} \tilde {\cal E}^D_\varphi\textmap{\varphi^D}{\cal O}_{S\times X}\to 0 
\end{equation}
(see Definitions \ref{iota-enlarged-def}, \ref{enlarged-ext}). Note that, in order to save on notations, we will use the superscript $D$ instead of $D_S$ in the notations ${\cal E}^D_\varphi$, $\varphi^D$, $\iota^D_\varphi$.

Taking into account Lemma \ref{R1=0}, we obtain the short exact sequence of sheaves on $S$:
\begin{equation}\label{ex-seq-on-S}
0\to p_{S*}({\cal E}_\varphi(D_S))\textmap{\iota^D_\varphi}p_{S*}(\tilde {\cal E}^D_\varphi)\textmap{\varphi^D} {\cal O}_S\to 0.
\end{equation}

A  local right inverse $\sigma:{\cal O}_{U}\to p_{S*}(\tilde {\cal E}^D_\varphi)_U$  of the epimorphism $\varphi^D$ in (\ref{ex-seq-on-S}) can be regarded as a fibrewise right splitting of (\ref{enl-ext-on-S}) over $U\times X$ via the identifications
$$
\Hom({\cal O}_{U},p_{U*}(\tilde {\cal E}^D_{\varphi_U}))=H^0(U,p_{U*}(\tilde {\cal E}^D_{\varphi_U}))=H^0(U\times X, \tilde {\cal E}^D_{\varphi_U})=\Hom({\cal O}_{U\times X},\tilde {\cal E}^D_{\varphi_U}).$$

Let ${\cal S}_\varphi^D$ be the sheaf of local right inverses of the sheaf epimorphism $\varphi^D$ in   (\ref{ex-seq-on-S}). The short exact sequence (\ref{ex-seq-on-S}) shows that this sheaf is naturally a $p_{S*}({\cal E}_\varphi(D_S))$-torsor. The sheaf $p_{S*}({\cal E}_\varphi(D_S))$ which acts on this torsor is locally isomorphic to ${\cal O}_S(W_D)$. Indeed, note that any  isomorphism
$$\big( p_X^{U\times X*}({\cal E})\textmap{f\simeq} {\cal E}_{\varphi_U}\big)\in {\cal T}_\varphi(U)$$
 defines  sheaf isomorphisms on $U\times X$:
\begin{equation}\label{f(DU)}
\begin{split}
f(D_U)&:p_X^{U\times X*}({\cal E}(D))=p_X^{U\times X*}({\cal E})(D_U)\map  {\cal E}_{\varphi_U}(D_U),
\\
f(D_U)_{D_U}&: p_X^{U\times X*}({\cal E}(D)_D)=p_X^{U\times X*}({\cal E})(D_U)_{D_U}\map    {\cal E}_{\varphi_U}(D_U)_{D_U},
\end{split}
\end{equation}
so also isomorphisms (which will be denoted by the same symbols)
\begin{equation}\label{p*f(D)}
\begin{split}
p_{U*}(p_X^{U\times X*}({\cal E}(D))={\cal O}_U(W_D)\map p_{U*}({\cal E}_{\varphi_U}(D_U))=p_{S*}({\cal E}_\varphi(D_S))_U,
\\
p_{U*}(p_X^{U\times X*}({\cal E}(D)_D)={\cal O}_U(V_D)\map   p_{U*}({\cal E}_{\varphi_U}(D_U)_{D_U}) =p_{S*}({\cal E}_\varphi(D_S)_{D_S})_U.
\end{split}
\end{equation}
For the first equalities in the above formulae we have used Remark \ref{pSpX}.  

We put
$$
\Theta_\varphi^D\edf {\cal S}_\varphi^D\times {\cal T}_\varphi, 
$$
i.e. for any open subspace $U\subset S$ we set $\Theta_\varphi^D(U)\edf {\cal S}_\varphi^D(U)\times{\cal T}_\varphi(U)$ and we endow this set with a right ${\cal O}(U,G_D)$-action by putting
$$
(\sigma,f)\circ (v,g)\edf  (\sigma+(\iota_{\varphi_U}^D\circ f(D_U))(v)\,,\,f\circ g)
$$
for any $(v,g)\in {\cal O}(U,W_D)\rtimes {\cal O}(U,G)={\cal O}(U,W_D\rtimes G)={\cal O}(U,G_D) $. 
%

We now define an ${\cal O}_S(G_D)$-equivariant sheaf morphism $\rho_\varphi^D:\Theta_\varphi^D\to {\cal O}_S(V_D)$ as follows:

For any open subspace $U\subset X$ the short exact sequence 
\begin{equation}\label{enlarged-ext-U}
0\to {\cal E}_{\varphi_U}(D_U)\textmap{\iota^D_{\varphi_U}} \tilde {\cal E}^D_{\varphi_U}\textmap{\varphi^D_U}{\cal O}_{U\times X}\to 0 
\end{equation}
is the $D_U$-enlarged extension of $0\to {\cal E}_{\varphi_U}\hookrightarrow \tilde {\cal E}_U\textmap{\varphi_U}{\cal O}_{U\times X}\to 0$, so it fits in the following commutative diagram with exact rows and columns:
\begin{equation}\label{restriction-diag-withS}
\begin{tikzcd}
&0\ar[d] &0 \ar[d] \\[-6pt]
0 \ar[r] & {\cal E}_{\varphi_U}  \ar[r, hook]  \ar[d, hook'] & \tilde {\cal E}_U   \ar[r, "\varphi_U" ]\ar[d, "j^D_{\varphi_U}"'] & {\cal O}_{U\times X} \ar[r]\ar[d, equal]& 0 \\
0  \ar[r]  & \mathcal{E}_{\varphi_U}(D_U)  \ar[d] \ar[r, "\iota^D_{\varphi_U}"] &  \tilde{\mathcal{E}}_{\varphi_U}^D \ar[d, "q_D^{\varphi_U}
"] \ar[r, "\varphi^D_U"'] & \mathcal{O}_{U\times X} \ar[l,"\sigma"', dashed ,bend right =40] \ar[r] & 0 \\
 & {\cal E}_{\varphi_U}(D_U)_{D_U} \ar[r, equal]\ar[d]&  {\cal E}_{\varphi_U}(D_U)_{D_U}  \ar[d] & & \\[-6pt]
 & 0 & 0
\end{tikzcd}.
\end{equation}
For any $(\sigma,f)\in {\cal S}_\varphi^D(U)\times{\cal T}_\varphi(U)$ we put : 
\begin{equation}\label{formula-for-rho}
\rho_\varphi^D(\sigma,f)\edf   f(D_U)_{D_U}^{-1}(q^{\varphi_U}_D\circ \sigma).
\end{equation}
This is an element of $H^0(U\times X,p_X^{U\times X*}({\cal E})(D_U)_{D_U})$, which (using the second isomorphism in (\ref{f(DU)}), (\ref{p*f(D)})) is identified with 
$$H^0(U\times X,p_X^{U\times X*}({\cal E}(D)_D))={\cal O}(U, H^0(X,{\cal E}(D)_D))={\cal O}(U,V_D).$$
The ${\cal O}_S(G_D)$-equivariance of $\rho_\varphi^D$ is straightforward.
\begin{dt} Our functor $\tensor[_S]{\Tg}{}$ is defined by $\tensor[_S]{\Tg}{}(\tilde{\cal E}\textmap{\varphi}{\cal O}_{S\times X})\edf (\Theta_\varphi^D,\rho_\varphi^D)$.	
\end{dt}

\subsection{The functor \texorpdfstring{$\tensor[_S]{\Ag}{}:\tensor[_S]{\Torss}{}^{\alpha_D}_{G_D}\to  \tensor[_S]{\Afff}{}_{E}$}{AgS} and its properties}
\label{functorAGS}
{\ }\\

Let $U\subset S$ be an open subspace and  let 
$$\beta\in H^0(U\times X\,,\,p_X^{U\times X*}({\cal E})(D_U)_{D_U})=H^0(U\,,\,p^{U\times X}_{U*}(p_X^{U\times X*}({\cal E}(D)_D))={\cal O}(U,V_D).$$
(see Remark \ref{pSpX}). We define the sheaf epimorphism 
$$u_\beta:p_X^{U\times X*}({\cal E})(D_U)\oplus {\cal O}_{U\times X}\to p_X^{U\times X*}({\cal E})(D_U)_{D_U}$$
by
$$
u_\beta(\eta,\zeta)\edf \eta_{D_U}+\zeta\beta,
$$
and we put $\tilde {\cal E}_\beta\edf \ker(u_\beta)$. In other words $\tilde {\cal E}_\beta$ is  the elementary transformation associated with $u_\beta$ (see \cite[1.2 p. 244]{Ma} for the general definition and properties of elementary transformations). We denote by $\varphi_\beta:\tilde {\cal E}_\beta\to {\cal O}_{U\times X}$ the epimorphism induced by the projection of $p_X^{U\times X*}({\cal E})(D_U)\oplus {\cal O}_{U\times X}$ on the second factor.  


For any  $\beta\in {\cal O}(U,V_D)$ we obtain the commutative diagram $(\Delta_\beta)$ with exact rows and columns 
\begin{equation}\label{Diag-beta-S}
\begin{tikzcd}[row sep=6mm, column sep=4mm]
&0\ar[d]&0\ar[d]&\\
0 \ar[r] & p_X^{U\times X*}({\cal E}) \ar[r ]  \ar[d, hook'] & \tilde{\cal E}_{\beta}  \ar[r, "\varphi_{\beta}" ]\ar[d,hook'] & {\cal O}_{U\times X} \ar[r]\ar[d, equal]& 0 \\
0  \ar[r]  & p_X^{U\times X*}({\cal E})(D_{U})\ar[d]   \ar[r, hook] &  p_X^{U\times X*}({\cal E})(D_{U})\oplus \mathcal{O}_{U\times X} \ar[d,"u_{\beta}"]   \ar[r, "p_2"] & \mathcal{O}_{U\times X}  \ar[r] & 0 \\
&p_X^{U\times X*}({\cal E})(D_{U})_{D_{U}}\ar[r, equal]\ar[d]& p_X^{U\times X*}({\cal E})(D_{U})_{D_{U}}\ar[d]\\
&0 & 0
 \end{tikzcd} \tag{$\Delta_\beta$}
 \end{equation}
in which the second row is canonically identified with the enlarged extension of the first, and, via this identification, $u_\beta$ coincides with the canonical epimorphism on $p_X^{U\times X*}({\cal E})(D_{U})\oplus \mathcal{O}_{U\times X}$ (see Definition \ref{enlarged-ext}). 
\vspace{2mm} 

We  need the   equivariance property of the correspondences $\beta\mapsto u_\beta$, $\beta\mapsto  \tilde {\cal E}_\beta$ with respect to the natural action of the group ${\cal O}(U,G_D)$. Let  $\gamma=(v,g)\in {\cal O}(U,G_D)$. We define an automorphism
$$
R_\gamma\in \Aut(p_X^{U\times X*}({\cal E})(D_U) \oplus {\cal O}_{U\times X})
$$
as follows: By the representability Theorem \ref{repres-Aut}, $g$ can be viewed as an automorphism of $p^{U\times X*}_X({\cal E})$. We will use the same symbol $g$ for this automorphism, and we will denote by $g(D_U)$, $g(D_U)_{D_U}$ the induced automorphisms of the sheaves
$
p_X^{U\times X*}({\cal E})(D_U)\, , \ p_X^{U\times X*}({\cal E})(D_U)_{D_U}$.
Similarly, by Remark \ref{pSpX}, $v\in {\cal O}(U,W_D)$ can be regarded as an element of 
\begin{equation*}
\begin{split}
H^0\big (U\times X, p^{U\times X*}_X({\cal E}(D))\big)&=H^0\big (U\times X, p_X^{U\times X*}({\cal E})(D_U)\big)\\
 &= \Hom\big({\cal O}_{U\times X},p_X^{U\times X*}({\cal E})(D_U)\big).
\end{split}
\end{equation*}
Therefore, the matrix 
$$
R_\gamma=\begin{pmatrix}
g(D_U) & v\\
0& 1	
\end{pmatrix}
$$
can be interpreted as an automorphism $R_\gamma$ of $p_X^{U\times X*}({\cal E})(D_U) \oplus {\cal O}_{U\times X}$. Direct computations show that:
%
%
\begin{pr}\label{equiv-u-beta-S}
For any $\gamma=(v,g)\in {\cal O}(U,G_D)$, $\beta\in {\cal O}(U,V_D)$ we have
$$
u_{\gamma \beta}\circ R_{\gamma}=g(D_U)_{D_U}\circ u_\beta.
$$
In other words, putting 
$$\gamma=(v,g),\ \gamma_0=g,\ R^0_\gamma=g(D_U)_{D_U}\in \Aut(p_X^{U\times X*}({\cal E})(D_U)_{D_U}),$$
we have the following commutative diagram: 
\begin{equation}\label{GD-equiv-S}
\begin{tikzcd}
p_X^{U\times X*}(\mathcal{E})(D_U) \oplus \mathcal{O}_{U\times X} \ar[r, "R_\gamma", "\simeq"'] \ar[d, "u_\beta"'] & p_X^{U\times X*}(\mathcal{E})(D_U)\oplus \mathcal{O}_{U\times X} \ar[d, "u_{\gamma \beta}"] \\
p_X^{U\times X*}(\mathcal{E})(D_U)_{D_U} \ar[r, "R_\gamma^0"', "\simeq"] & p_X^{U\times X*}(\mathcal{E})(D_U)_{D_U}	
\end{tikzcd}
\end{equation}
with isomorphic rows. 

In particular $R_\gamma$ maps isomorphically $\tilde {\cal E}_\beta\edf \ker(u_\beta)$ onto $\tilde {\cal E}_{\gamma\beta}\edf \ker(u_{\gamma\beta})$.
\end{pr}

Let now $(\Theta,\rho)$ be an object of the groupoid $\tensor[_S]{\Torss}{}^{\alpha_D}_{G_D}$. For a local section $\theta \in \Theta(U_\theta)$ we obtain the locally free sheaf $\tilde {\cal E}_{\rho(\theta)}$ on $U_\theta\times X$. 

For two local sections $\theta\in \Theta(U_\theta)$, $\theta'\in \Theta(U_{\theta'})$ of $\Theta$ we define $\gamma^{\theta'}_\theta\in {\cal O}(U_\theta\cap U_{\theta'},G_D)$ by the formula
\begin{equation}\label{gamma-theta-theta'}
\theta|_{U_\theta\cap U_{\theta'}}=\theta'|_{U_\theta\cap U_{\theta'}}\,  \gamma^{\theta'}_\theta.
\end{equation}

This implies 
$$\rho(\theta')|_{U_\theta\cap U_{\theta'}}=\gamma^{\theta'}_\theta\rho(\theta)|_{U_\theta\cap U_{\theta'}},$$
so, by Proposition \ref{equiv-u-beta-S}, it follows that $R_{\gamma^{\theta'}_\theta}$ maps isomorphically $\tilde {\cal E}_{\rho(\theta)}|_{(U_\theta\cap U_{\theta'})\times X}$ onto $\tilde {\cal E}_{\rho(\theta')}|_{(U_\theta\cap U_{\theta'})\times X}$, so it induces an isomorphism
\begin{equation}\label{phi-theta-theta'}
\varphi^{\theta'}_\theta:\tilde {\cal E}_{\rho(\theta)}|_{(U_\theta\cap U_{\theta'})\times X}\textmap{\simeq}\tilde {\cal E}_{\rho(\theta')}|_{(U_\theta\cap U_{\theta'})\times X}.
\end{equation}

Applying the gluing functor given by Proposition \ref{Gluing} to the open cover 
$$(U_\theta\times X)_{\substack{\theta\ {\rm loc.\ sect.}\\{\rm  \ of}\ \Theta}}$$
of $S\times X$ and the pair $\big ((\tilde {\cal E}_{\rho(\theta)})\,,\, (\varphi^{\theta'}_\theta)\big)$, we obtain a locally free sheaf $\tilde {\cal E}_{\Theta,\rho}$ on $S\times X$.

By the universal property given by Proposition \ref{Gluing}, we see that the family of sheaf epimorphisms 
$$
\begin{tikzcd}
\tilde {\cal E}_{\rho(\theta)}\ar[r,hook]\ar[rr, bend right=15, "\varphi_{\rho(\theta)}"']& p_X^*({\cal E})(D_{U_\theta})\oplus {\cal O}_{U_\theta\times X}\ar[r, two heads]& {\cal O}_{U_\theta\times X}.
\end{tikzcd}
$$
defines canonical epimorphism $\varphi_{\Theta,\rho}:\tilde {\cal E}_{\Theta,\rho}\to {\cal O}_{U_\theta\times X}$.
\begin{dt}\label{def-AgS}
Our functor $\tensor[_S]{\Ag}{}$ is defined by 
$\tensor[_S]{\Ag}{}(\Theta,\rho)\edf \big(\tilde {\cal E}_{\Theta,\rho}\textmap{\varphi_{\Theta,\rho}}{\cal O}_{S\times X}\big)$.
\end{dt}
For the construction of a natural isomorphism $\id\textmap{\simeq}\tensor[_S]{\Tg}{}\circ \tensor[_S]{\Ag}{}$ we will need
%
%
an explicit description of the following objects: 
\begin{enumerate}
\item 	The kernel ${\cal E}_{\varphi_{\Theta,\rho}}$ of the augmented bundle $\varphi_{\Theta,\rho}:\tilde {\cal E}_{\Theta,\rho}\to {\cal O}_{S\times X}$.
\item The $D_S$-enlarged extension 
$$
0\to  {\cal E}_{\varphi_{\Theta,\rho}}(D_S) \textmap{\iota^D_{\varphi_{\Theta,\rho}}}   \tilde {\cal E}^D_{\varphi_{\Theta,\rho}}\to {\cal O}_{S\times X}\to 0
$$
of  the extension  $
0\to {\cal E}_{\varphi_{\Theta,\rho}}\hookrightarrow \tilde {\cal E}_{\Theta,\rho}\textmap{\varphi_{\Theta,\rho}} {\cal O}_{S\times X}\to 0$.
\item The canonical epimorphism 
$$q_D^{\varphi_{\Theta,\rho}}: \tilde {\cal E}^D_{\varphi_{\Theta,\rho}}\to {\cal E}_{\varphi_{\Theta,\rho}}(D_S)_{D_S}$$ 
on the central term of this enlarged extension. 
\end{enumerate}
%


We will see that these three  objects are  associated with gluing data parameterised by local sections $\theta$ of $\Theta$.\\

Let  $\theta$ be a local section  of $\Theta$.  Taking $\beta=\rho(\theta)$ in  the diagram (\ref{Diag-beta-S}) of this section,  we obtain the commutative diagram $(\Delta_\theta)\edf (\Delta_{\rho(\theta)})$.
Let $\theta$, $\theta'$ be local sections  of $\Theta$. Put $U_{\theta\theta'}\edf U_\theta\cap U_{\theta'}$, $D_\theta\edf D_{U_\theta}$, $D_{\theta\theta'}\edf D_{U_{\theta\theta'}}$ and  denote by   $p^{\theta}_X$, $p^{\theta\theta'}_X$ the projections $p^{{U_\theta}\times X}_X:{U_\theta}\times X\to X$, $p^{U_{\theta\theta'}\times X}_X:{U_{\theta\theta'}}\times X\to X$ respectively to save on notations. Let $\gamma^{\theta'}_\theta=(v^{\theta'}_\theta,g^{\theta'}_\theta)\in {\cal O}(U_{\theta\theta'},G_D)$ be given by (\ref{gamma-theta-theta'}).
\begin{dt}\label{DefRg-gamma} 
We define the isomorphism of commutative diagrams (in the obvious sense)
$$\Rg_{\gamma^{\theta'}_\theta}: (\Delta_{\theta})|_{U_{\theta\theta'}\times X}\to (\Delta_{\theta'})|_{U_{\theta\theta'}\times X}$$
as follows:
\begin{enumerate}
\item On the sheaves $p_X^{\theta\theta'*}({\cal E})$, $p_X^{\theta\theta'*}({\cal E})(D_{\theta\theta'})$, $p_X^{\theta\theta'*}({\cal E})(D_{\theta\theta'})_{D_{\theta\theta'}}$ the isomorphism	 $\Rg_{\gamma^{\theta'}_\theta}$ acts by $g^{\theta'}_\theta$,   $g^{\theta'}_\theta(D_{\theta\theta'})$, respectively $g^{\theta'}_\theta(D_{U_{\theta\theta'}})_{D_{\theta\theta'}}$.
\item On the sheaves  $p_X^{\theta\theta'*}({\cal E})(D_{\theta\theta'})\oplus {\cal O}_{X\times U_{\theta\theta'}}$, $\tilde {\cal E}_{\rho(\theta)}|_{U_{\theta\theta'}\times X}$ the isomorphism	 $\Rg_{\gamma^{\theta'}_\theta}$ acts by $R_{\gamma^{\theta'}_\theta}$, respectively by its restriction to $\tilde {\cal E}_{\rho(\theta)}|_{U_{\theta\theta'}\times X}$.
\item On ${\cal O}_X$ the isomorphism	$\Rg_{\gamma^{\theta'}_\theta}$ acts by identity.
\end{enumerate}
\end{dt}

Using Remark \ref{FunctorialityGluing} we obtain a commutative diagram with exact rows and columns
{\small 
\begin{equation*}
\begin{tikzcd}[row sep=8mm, column sep=3mm]
& 0\ar[d] & 0\ar[d]\\
0 \ar[r] & \hspace*{-4mm} {\displaystyle\coprod_{\substack{\theta{\rm\, loc.\, sect.}\\ {\rm of}\,\Theta}}}\hspace*{-1mm} p^{\theta*}_X({\cal E})\big/(\Rg_{\gamma^{\theta'}_{\theta}}) \ar[r, ]  \ar[d, ] & \tilde{\cal E}_{\Theta,\rho}=\hspace*{-3mm} {\displaystyle\coprod_{\substack{\theta{\rm\, loc.\, sect.}\\ {\rm of}\,\Theta}}}\hspace*{-1mm}\tilde{\cal E}_{\rho(\theta)}\big/(\Rg_{\gamma^{\theta'}_{\theta}})  \ar[r, "\varphi_{\Theta,\rho}" ]\ar[d,hook'] & {\cal O}_{S\times X} \ar[r]\ar[d, equal]& 0 \\
0  \ar[r]  &\hspace*{-4mm} {\displaystyle\coprod_{\substack{\theta{\rm\, loc.\, sect.}\\ {\rm of}\,\Theta}}}\hspace*{-1mm} p^{\theta*}_X({\cal E})(D_{\theta})\big/(\Rg_{\gamma^{\theta'}_{\theta}})   \ar[r, "\jg"] \ar[d] & \hspace*{-4mm} {\displaystyle\coprod_{\substack{\theta{\rm\, loc.\, sect.}\\ {\rm of}\,\Theta}}}\hspace*{-1mm}p^{\theta*}_X({\cal E})(D_{\theta})\oplus \mathcal{O}_{U_\theta\times X})\big/(\Rg_{\gamma^{\theta'}_{\theta}})   \ar[r, "\pg "] \ar[d,"{\displaystyle\coprod_{\substack{\theta{\rm\, loc.\, sect.}\\ {\rm of}\,\Theta}}}\hspace*{-1mm} u_{\rho(\theta)}"] & \mathcal{O}_{S\times X}  \ar[r] & 0  \\
&{\displaystyle\coprod_{\substack{\theta{\rm\, loc.\, sect.}\\ {\rm of}\,\Theta}}}\hspace*{-1mm} p^{\theta*}_X({\cal E})(D_{\theta})_{D_{\theta}}\big/(\Rg_{\gamma^{\theta'}_{\theta}})\ar[r, equal]\ar[d]&{\displaystyle\coprod_{\substack{\theta{\rm\, loc.\, sect.}\\ {\rm of}\,\Theta}}}\hspace*{-1mm} p^{\theta*}_X({\cal E})(D_{\theta})_{D_{\theta}}\big/(\Rg_{\gamma^{\theta'}_{\theta}})\ar[d]&&\\
&0&0&
 \end{tikzcd}  
\end{equation*} }
of sheaves associated to gluing data, with canonical isomorphisms
\begin{equation*}
\begin{array}{l}
{\displaystyle\coprod_{\substack{\theta{\rm\, loc.\, sect.}\\ {\rm of}\,\Theta}}}p^{\theta*}_X({\cal E})\big/(\Rg_{\gamma^{\theta'}_{\theta}}) =\ker(\varphi_{\Theta,\rho})\eqcolon{\cal E}_{\varphi_{\Theta,\rho}}, \\ {\displaystyle\coprod_{\substack{\theta{\rm\, loc.\, sect.}\\ {\rm of}\,\Theta}}} p^{\theta*}_X({\cal E})(D_{\theta})\big/(\Rg_{\gamma^{\theta'}_{\theta}}) ={\displaystyle\coprod_{\substack{\theta{\rm\, loc.\, sect.}\\ {\rm of}\,\Theta}}}p^{\theta*}_X({\cal E})\big/(\Rg_{\gamma^{\theta'}_{\theta}})(D_S)={\cal E}_{\varphi_{\Theta,\rho}}(D_S).
\end{array}
\end{equation*}
Note that via the above identifications, the upper vertical arrow on the left is the canonical monomorphism ${\cal E}_{\varphi_{\Theta,\rho}}\hookrightarrow{\cal E}_{\varphi_{\Theta,\rho}}(D_S)$.

  Put  $\tilde {\cal E}\edf \tilde {\cal E}_{\Theta,\rho}$,  $\varphi\edf \varphi_{\Theta,\rho}$ to save on notations. Using the universal property of the enlarged extension (see Remark \ref{enlarged-remark1}) and the characterisation of the canonical epimorphism (see Remark \ref{PropOfCanonEpim}), we obtain the following proposition which will be  used for the construction of the natural isomorphism $\id_{\Torsss_{G_D}^{\alpha_D}}\stackrel{\simeq}{\to}  \tensor[_S]{\Tg}{}\circ\tensor[_S]{\Ag}{}$:
%


%

\begin{pr} \label{propertiesAg-S}  Let $(\Theta,\rho)$ be an object of the groupoid $\Torsss_{G_D}^{\alpha_D}$.
 Put  $\tilde {\cal E}\edf \tilde {\cal E}_{\Theta,\rho}$,  $\varphi\edf \varphi_{\Theta,\rho}$ to save on notations.
With the notations and under assumptions above, we have:
\begin{enumerate}
	\item The kernel 	${\cal E}_{\varphi}$ of the augmented bundle $\tilde {\cal E}\textmap{\varphi} {\cal O}_{S\times X}$ defined above is canonically identified with the sheaf ${\displaystyle\coprod_{\substack{\theta{\rm\, loc.\, sect.}\\ {\rm of}\,\Theta}}}p^{\theta*}_X({\cal E})\big/(\Rg_{\gamma^{\theta'}_{\theta}})$ associated with  the ${\cal O}_{S\times X}$-module gluing data 
$\big(({p^{\theta*}_X({\cal E}))_{\substack{\theta{\rm\, loc.\, sect.}\\ {\rm   of}\,\Theta}}}\,,\,(\Rg_{\gamma^{\theta'}_{\theta}})\big)$.
\item The $D_S$-enlarged extension 
$$0\to  {\cal E}_{\varphi}(D_S) \textmap{\iota^D_\varphi}   \tilde {\cal E}^D_{\varphi}\to {\cal O}_{S\times X}\to 0
$$
of the extension
$$
0\to {\cal E}_{\varphi}\hookrightarrow \tilde {\cal E}\textmap{\varphi}{\cal O}_{S\times X}\to 0
$$
is canonically identified with the extension
  {\small 
\begin{equation*}
\begin{tikzcd}[row sep=9mm, column sep=4mm]
0  \ar[r]  &\hspace*{-4mm} {\displaystyle\coprod_{\substack{\theta{\rm\, loc.\, sect.}\\ {\rm of}\,\Theta}}} p^{\theta*}_X({\cal E})(D_{\theta})\big/(\Rg_{\gamma^{\theta'}_{\theta}})   \ar[r, "\jg"]  & \hspace*{-4mm} {\displaystyle\coprod_{\substack{\theta{\rm\, loc.\, sect.}\\ {\rm of}\,\Theta}}}p^{\theta*}_X({\cal E})(D_{\theta})\oplus \mathcal{O}_{U_\theta\times X})\big/(\Rg_{\gamma^{\theta'}_{\theta}})   \ar[r, "\pg "]  & \mathcal{O}_{S\times X}  \ar[r] & 0   \end{tikzcd}  
\end{equation*} }
 of sheaves associated with gluing data, where the monomorphism $\jg$ is obtained by gluing (in the sense of Remark \ref{Mor-Ass-Mor-GluingData})   the  family of inclusions 
 $$
 \underset{\substack{{\rm\, loc.\, sect.}\\ {\rm   of}\,\Theta}}{\theta}  \mapsto (j_1:p_X^{\theta*}({\cal E})(D_{U_\theta})\hookrightarrow p_X^{\theta*}({\cal E})(D_{U_\theta})\oplus {\cal O}_{U_\theta\times X}),
 $$
 and the  epimorphism $\pg$ is obtained by gluing (in the sense of Remark \ref{Mor-Ass-Mor-GluingData})  the family  of projections 
 $$\underset{\substack{{\rm\, loc.\, sect.}\\ {\rm   of}\,\Theta}}{\theta}\mapsto (p_2:p_X^{\theta*}({\cal E})(D_{U_\theta})\oplus {\cal O}_{U_\theta\times X}\to {\cal O}_{U_\theta\times X}).$$
\item Via the identification $\tilde {\cal E}= {\displaystyle\coprod_{\substack{\theta{\rm\, loc.\, sect.}\\ {\rm of}\,\Theta}}(p_X^{\theta*}({\cal E})(D_{U_\theta})\oplus \mathcal{O}_{U_\theta\times X})}\big/(\Rg_{\gamma^{\theta'}_{\theta}})$\,, the canonical epimorphism 
$q_D^\varphi:\tilde {\cal E}_\varphi^D\to {\cal E}_\varphi(D_S)_{D_S}$
coincides with the morphism ${\displaystyle\coprod_{\substack{\theta{\rm\, loc.\, sect.}\\ {\rm of}\,\Theta}}} u_{\rho(\theta)}$ obtained by gluing the family of morphisms 
$$\big(  p_X^{\theta*}({\cal E})(D_{U_\theta})\oplus \mathcal{O}_{U_\theta\times X} \textmap{u_{\rho(\theta)}} p_X^{\theta*}({\cal E})(D_{U_\theta})_{D_{U_\theta}}\big)_{\substack{\theta{\rm\, loc.\, sect.}\\ {\rm of}\,\Theta}} $$
\end{enumerate}
\end{pr}

Using again the properties of the sheaf associated with gluing data stated and the morphisms associated with gluing data morphisms (Proposition \ref{Gluing}, Remarks \ref{Mor-Ass-Mor-GluingData}, \ref{Glue-sheaves-Ox-modules}, \ref{FunctorialityGluing}), we obtain:  
\begin{enumerate}
	\item  A family 
	$$(\psi_\theta:p_X^{\theta*}({\cal E})\textmap{\simeq} {\cal E}_{\varphi_{U_\theta}}={\cal E}_\varphi|_{U_\theta\times X})_{\substack{\theta{\rm\, loc.\, sect.}\\ {\rm of}\,\Theta}}$$
of canonical sheaf isomorphisms  satisfying the compatibility condition
\begin{equation}\label{Compat-psi-S}
	 \psi_{\theta'}|_{U_{\theta\theta'}\times X}=\psi_\theta|_{U_{\theta\theta'}\times X}\circ \Rg_{\gamma^{\theta'}_\theta} \hbox{ for any local sections }\theta,\ \theta' \hbox{ of }\Theta. 
\end{equation}
	\item A family   
$$(\Psi_\theta:p_X^{\theta*}({\cal E})(D_{\theta})\oplus {\cal O}_{U_\theta\times X}\textmap{\simeq}\tilde {\cal E}^D_{\varphi_{U_\theta}}=\tilde{\cal E}_\varphi^D|_{U_\theta\times X})_{\substack{\theta{\rm\, loc.\, sect.}\\ {\rm of}\,\Theta}}$$
of canonical sheaf isomorphisms making commutative the diagrams
\begin{equation}\label{comm-diag-Psi-theta1-S}
\begin{tikzcd}[column sep=4mm]
 0\ar[r]&p_X^{\theta*}({\cal E})(D_{\theta})\ar[r,"j_1"]\ar[d,"\psi_\theta(D)"]&p_X^{\theta*}({\cal E})(D_{\theta})\oplus {\cal O}_{U_\theta\times X}\ar[r,"p_2"]\ar[d, "\Psi_\theta"]&{\cal O}_{U_\theta\times X}\ar[r]\ar[d,equal]&0\\
 0\ar[r]&{\cal E}_{\varphi_{U_\theta}}(D_{\theta})\ar[r, "\iota_{\varphi_{U\theta}}^D"]& \tilde {\cal E}_{\varphi_{U_\theta}}^D\ar[r,"\varphi_{U\theta}"]&{\cal O}_{U_\theta\times X}\ar[r] &0
\end{tikzcd},\end{equation}
\begin{equation}\label{comm-diag-Psi-theta2-S} 
\begin{tikzcd}[row sep=6mm, column sep=18mm]
p_X^{\theta*}({\cal E})(D_{\theta})\oplus {\cal O}_{U_\theta\times X}\ar[d, "u_{\rho(\theta)}"']\ar[r, "\Psi_\theta"] & \tilde {\cal E}^D_{\varphi_{U_\theta}}\ar[d, "q^\varphi_D"]\\	
p_X^{\theta*}({\cal E})(D_{\theta})_{D_{\theta}} \ar[r, "\psi_\theta(D_{\theta})_{D_{\theta}}"]&{\cal E}_{\varphi_{U_\theta}}(D_{\theta})_{D_{\theta}}
\end{tikzcd}
\end{equation}		
and  satisfying the compatibility condition
\begin{equation}\label{Compat-Psi-S}
\Psi_{\theta}|_{(U_\theta\cap U_{\theta'})\times X}{=\Psi_{\theta'}|_{(U_\theta\cap U_{\theta'})\times X}\circ R_{\gamma^{\theta'}_\theta}} \hbox{ for any local sections }\theta,\ \theta' \hbox{ of }\Theta.
\end{equation}
\end{enumerate}

\subsection{The natural isomorphisms \texorpdfstring{$\tensor[_S]{\Ag}{}\circ\tensor[_S]{\Tg}{}\to \id$}{agtgs}, \texorpdfstring{$\id\to \tensor[_S]{\Tg}{}\circ\tensor[_S]{\Ag}{}$}{tgsags}} 
\subsubsection{The natural isomorphism \texorpdfstring{$\tensor[_S]{\Ag}{}\circ\tensor[_S]{\Tg}{}\to \id$}{agtgs-sub}}
{\ }
For any object $\tilde{\cal E}\textmap{\varphi}{\cal O}_{S\times X}$ of $\tensor[_S]{\Ag}{}\circ\tensor[_S]{\Tg}{}$ we will construct a {\it canonical} isomorphism 
$$
\tensor[_S]{\Ag}{}(\tensor[_S]{\Tg}{}(\tilde{\cal E}\textmap{\varphi}{\cal O}_X))\textmap{\simeq}(\tilde{\cal E}\textmap{\varphi}{\cal O}_{S\times X})
$$
of augmented bundles. Let $\tilde {\cal E}\textmap{\varphi}{\cal O}_{S\times X}$ be an object of the groupoid $\Afffs_E$. Recall that $\tensor[_S]{\Tg}{}(\tilde{\cal E}\textmap{\varphi}{\cal O}_X)=(\Theta_\varphi^D,\rho _\varphi^D)$ (see section \ref{functorTGS}).

Put $\Theta\edf \Theta_\varphi^D$, $\rho\edf \rho _\varphi^D$ to save on notations.
The augmented bundle  $\tensor[_S]{\Ag}{}(\Theta,\rho)$ on $S\times X$ is defined using the gluing principle for sheaves   by: 
$$
\tensor[_S]{\Ag}{}(\Theta,\rho)= \bigg(\coprod_{\substack{\theta{\rm\, loc.\, sect.}\\ {\rm of}\,\Theta}} \tilde{{\cal E}}_{\rho(\theta)} /_{(\varphi_\theta^{\theta'})}\textmap{\varphi_{\Theta,\rho}} \mathcal{O}_X\bigg).
$$
By the universal property given by Proposition \ref{Gluing} it suffices  to construct, for any local section $\theta$ of $\Theta$, an isomorphism of augmented bundles 
\begin{equation}\label{isom-of-augmented-S}
\begin{tikzcd}[row sep=5mm]
\tilde{\cal E}_{\rho(\theta)} \ar[dr, "\varphi_{\rho(\theta)}"']\ar[rr, "e_\theta"]& & \tilde {\cal E}_{U_\theta} \ar[dl, "\varphi_{U_\theta}"]	\\
& \mathcal{O}_{U_\theta\times X}
\end{tikzcd}
\end{equation}
such that the following compatibility condition is satisfied
\begin{equation}\label{comp-for-e-beta-S}
e_{\theta'}|_{(U_\theta\cap U_{\theta'})\times X}\circ \varphi_\theta^{\theta'}=e_\theta|_{(U_\theta\cap U_{\theta'})\times X} \hbox{ for local sections }\theta,\ \theta' \hbox{ of } \Theta.
\end{equation}
%

The idea of the construction of $e_\theta$ is very natural: A local section $\theta$ of $\Theta=\Theta_\varphi^D$ is a pair $(\sigma,f)$ consisting of a splitting $\sigma $ of the exact sequence 
$$
0 \to \mathcal{E}_{\varphi_{U_\theta}}(D_\theta) \textmap{\iota^D_\varphi} \tilde{\mathcal{E}}_{\varphi_{U_\theta}}^D\textmap{\varphi^D} \mathcal{O}_{U_\theta\times X}\to 0. \eqno{(e_{\varphi_{U_\theta}}^D)}
$$
and an isomorphism of ${\cal O}_{U_\theta\times X}$-modules  $p_X^{\theta*}({\cal E})\textmap{f\simeq}{\cal E}_{\varphi_{U_\theta}}$.

 We will use $\sigma$ to identify $\tilde{\mathcal{E}}_{\varphi_{U_\theta}}^D$ with $\mathcal{E}_{\varphi_{U_\theta}}(D_\theta)\oplus \mathcal{O}_{U_\theta\times X}$ and $f\oplus\id$ to further identify this direct sum to $p_X^{\theta*}({\cal E})(D_\theta)\oplus \mathcal{O}_{U_\theta\times X}$. This will allow us to identify the subsheaf $\tilde {\cal E}_{U_\theta}$ of $\tilde{\mathcal{E}}_{\varphi_{U_\theta}}^D$ to precisely the elementary transformation 
 $$\tilde {\cal E}_{\rho(\theta)}\edf \ker(u_{\rho(\theta)})\hookrightarrow p_X^{\theta*}({\cal E})(D_\theta)\oplus \mathcal{O}_{U_\theta\times X}.$$

More precisely, let $\theta\edf(\sigma, f)$ be as above and $\beta \edf \rho(\theta)=f(D_\theta)_{D_\theta}^{-1}(q_D^{\varphi_{U_\theta}}\circ \sigma )$ (see formulae (\ref{f(DU)})), (\ref{formula-for-rho}) and the diagram (\ref{restriction-diag-withS})). The   diagram 
\begin{equation}\label{lower-diagr}
\begin{tikzcd}[column sep =2cm]
p_X^{\theta*}({\cal E})(D_\theta)\oplus \mathcal{O}_{U_\theta\times X}  \ar[d, "u_\beta"'] \ar[r, "(\iota_{\varphi_{U_\theta}}^D \circ f(D_\theta){\,,\,} \sigma)", "\simeq"'] &   \tilde{\mathcal{E}}_{\varphi_{U_\theta}}^D \ar[d, "q_D^{\varphi_{U_\theta}}"]      \\
p_X^{\theta*}({\cal E})(D_\theta)_{D_\theta} \ar[r, "f(D_\theta)_{D_\theta}", "\simeq"'] \ar[d]&    {\cal E}_{\varphi_{U_\theta}}(D_\theta)_{D_\theta}\ar[d]   \\
 0 &   0
\end{tikzcd}
\end{equation}
is commutative. Indeed, we have
$$
f(D_\theta)_{D_\theta}\circ u_\beta (s, \zeta)=f(D_\theta)_{D_\theta}(s_{D_\theta}+\zeta \beta)=(f(D_\theta)(s))_{D_\theta}+\zeta q_D^{\varphi_{U_\theta}} \circ \sigma,
$$
and, using the commutativity of (\ref{restriction-diag}), we also obtain 
\begin{equation*}
\begin{split}
q_D^{\varphi_{U_\theta}}  \circ (\iota_{\varphi_{U_\theta}}^D \circ f(D_\theta),\sigma)(s,\zeta)&= q_D^{\varphi_{U_\theta}}\big(\iota_{\varphi_{U_\theta}}^D \circ f(D_\theta)(s)+\zeta \sigma \big)=\\
&=q_D^{\varphi_{U_\theta}} \big(\iota_{\varphi_{U_\theta}}^D \circ f(D_\theta)(s)\big)+ \zeta q_D^\varphi \circ \sigma=\\
& =(f(D_\theta)(s))_{D_\theta}+ \zeta q_D^{\varphi_{U_\theta}} \circ \sigma. 
\end{split}
\end{equation*}
It follows that the sheaf isomorphism
$$
I_\theta\edf (\iota_{\varphi_{U_\theta}}^D \circ f(D_\theta){\,,\,} \sigma): p_X^{\theta*}({\cal E})(D_\theta)\oplus \mathcal{O}_{U_\theta\times X}\textmap{\simeq} \tilde {\cal E}^D_{\varphi_{U_\theta}}
$$
maps isomorphically $\tilde {\cal E}_\beta\edf \ker(u_\beta)$ onto $\ker(q_D^{\varphi_{U_\theta}})=\tilde {\cal E}_{U_\theta}$, so it induces an isomorphism 
$$
e_\theta: \tilde {\cal E}_{\rho(\theta)}\textmap{\simeq}\tilde {\cal E}_{U_\theta}.
$$
In other words, the commutative diagram (\ref{lower-diagr}) can be extended  to the commutative diagram
\begin{equation}\label{big-diagr}
\begin{tikzcd}[row sep=6mm, column sep =1cm]
0 \ar[d] &   0\ar[d]\\
\tilde{\cal E}_\beta \ar[d, hook'] \ar[r, "e_\theta", "\simeq"']&   \tilde{\cal E}_{U_\theta}\ar[d, "j_{\varphi_{U_\theta}}^D"]\\
p_X^{\theta*}({\cal E})(D_\theta)\oplus \mathcal{O}_{U_\theta\times X}  \ar[d, "u_\beta"'] \ar[r, "I_\theta", "\simeq"'] &   \tilde{\cal E}_{\varphi_{U_\theta}}^D \ar[d, "q_D^{\varphi_{U_\theta}}"]    \\
p_X^{\theta*}({\cal E})(D_\theta)_{D_\theta}  \ar[r, "f(D_\theta)_{D_\theta}", "\simeq"'] \ar[d]&   {\cal E}_{\varphi_{U_\theta}}(D_\theta)_{D_\theta}\ar[d]   \\
 0 &    0
\end{tikzcd}
\end{equation}


 Now note that the diagram
\begin{equation}\label{diag-with-rectangle-triangles}
\begin{tikzcd}[row sep=4mm, column sep =1cm]
0 \ar[d] & &  0\ar[d]\\
\tilde{\cal E}_\beta \ar[dd, hook'] \ar[dr, "\varphi_\beta"'] \ar[rr, "e_\theta", "\simeq"']&  & \tilde{\cal E}_{U_\theta}\ar[dl, "\varphi_{U_\theta}"']\ar[dd, "j_{\varphi_{U_\theta}}^D"]\\
&{\cal O}_{U_\theta\times X}&\\
p_X^{\theta*}({\cal E})(D_\theta)\oplus \mathcal{O}_{U_\theta\times X}    \ar[rr, "I_\theta", "\simeq"'] \ar[ur, "p_2"] & &  \tilde{\mathcal{E}}_{\varphi_{U_\theta}}^D   \ar[ul, "\varphi^D_{U_\theta}"']   
\end{tikzcd}
\end{equation}
is commutative. Indeed, the commutativity of the lower triangle follows taking into account that $\varphi^D_{U_\theta}\circ \iota^D_{\varphi_{U_\theta}}=0$ and, since $\sigma$ is a right splitting of the enlarged extension $(e_{\varphi_{U_\theta}}^D)$, we have  $\varphi^D_{U_\theta}\circ\sigma=\id_{{\cal O}_{U_\theta\times X}}$.
Therefore, the commutativity of the upper triangle follows from the commutativity of the lateral  triangles and the commutativity of the rectangle. The latter is obvious since  $e_\theta$ is induced by $I_\theta$. The commutativity of the lateral triangles follows using the commutative diagrams (\ref{restriction-diag}), ($\Delta_\beta$).

The commutativity of the upper triangle yields $\varphi_{U_\theta}\circ e_\theta=\varphi_{\rho(\theta)}$, so $e_\theta$ is an isomorphism of augmented bundles as claimed in (\ref{isom-of-augmented-S}). 

It remains to prove the compatibility condition (\ref{comp-for-e-beta-S}). Let $\theta$, $\theta'$ be local sections of $\Theta$ and let again $\gamma^{\theta'}_\theta\in {\cal O}(U_{\theta\theta'},G_D)$ be defined by the condition
$\theta|_{U_{\theta\theta'}}=\theta'|_{U_{\theta\theta'}}\,  \gamma^{\theta'}_\theta$. Here we used again the simplifying  notations introduced at the end of section \ref{functorAGS}
 Now note that  in the following diagram
\begin{equation}\label{diagram-with-R-on-the-left}
\begin{tikzcd}[row sep=5mm, column sep =2.5cm]
p_X^{\theta*}({\cal E})(D_{\theta\theta'})\oplus \mathcal{O}_{U_{\theta\theta'}\times X}\ar[dd,"R_{\gamma^{\theta'}_\theta}"]\ar[dr, "I_\theta|_{U_{\theta\theta'}}" description]&&\tilde {\cal E}_{\rho(\theta)}\ar[dd, "\varphi^{\theta'}_\theta"] \ar[dl, "j^D_{\varphi_{U_{\theta\theta'}}}\circ e_\theta|_{U_{\theta\theta'}}" description]\ar[ll, hook']\\
& \tilde {\cal E}^D_{\varphi_{U_{\theta\theta'}}} &\\	
p_X^{\theta*}({\cal E})(D_{\theta\theta'})\oplus \mathcal{O}_{U_{\theta\theta'}\times X}\ar[ur, "I_{\theta'}|_{U_{\theta\theta'}}" description]&&\tilde {\cal E}_{\rho(\theta')}\ar[ul, "j^D_{\varphi_{U_{\theta\theta'}}}\circ e_{\theta'}|_{U_{\theta\theta'}}" description]\ar[ll,hook']
\end{tikzcd}
\end{equation}
the rectangle is commutative by the definition of  $\varphi^{\theta'}_\theta$, the upper and lower triangles are commutative by the definitions of $e_\theta$, $e_{\theta'}$. Therefore the commutativity of the right-hand triangle follows from the commutativity of the left-hand triangle. But,  putting $\theta=(\sigma,f)$,  $\theta'=(\sigma',f')$,  $\gamma^{\theta'}_\theta=(v^{\theta'}_\theta,g^{\theta'}_\theta)\in {\cal O}(U_{\theta\theta'},G_D)$, we have on $U_{\theta\theta'}\times X$:
$$
(\sigma,f)=(\sigma',f')(v^{\theta'}_\theta,g^{\theta'}_\theta)=\big(\sigma'+(\iota^D_{\varphi_{U_{\theta\theta'}}}\circ f'(D_{\theta\theta'}))(v)\,,\,f'\circ g^{\theta'}_\theta\big),
$$
so
\begin{equation*}
\begin{split}
I_{\theta'}|_{U_{\theta\theta'}}\circ R_{\gamma^{\theta'}_\theta}&=(\iota^D_{\varphi_{U_{\theta\theta'}}}\circ f'(D_{\theta\theta'})\,,\,\sigma')\begin{pmatrix}
g^{\theta'}_\theta(D_{\theta\theta'}) & v\\ 
0& 1	
\end{pmatrix}=(\iota^D_{\varphi_{U_{\theta\theta'}}}\circ f(D_{\theta\theta'})\,,\ \sigma)\\
&= I_\theta|_{U_{\theta\theta'}},
\end{split}
\end{equation*}
which proves that the left hand triangle is commutative, as claimed.
Therefore in conclusion  we have:

\begin{pr}
The family of augmented bundle isomorphisms 	
$$
\begin{tikzcd}[row sep=5mm]
\tilde{\cal E}_{\rho(\theta)} \ar[dr, "\varphi_{\rho(\theta)}"']\ar[rr, "e_\theta"]& & \tilde {\cal E}_{U_\theta} \ar[dl, "\varphi"]	\\
& \mathcal{O}_{U_\theta\times X}
\end{tikzcd}
$$
defined above verifies the compatibility condition (\ref{comp-for-e-beta-S}), so, by the universal property given by Proposition \ref{Gluing}, defines a canonical isomorphism
$$
\Ag(\Tg(\tilde{\cal E}\textmap{\varphi}{\cal O}_{S\times X}))\edf\bigg( {\displaystyle\coprod_{\substack{\theta{\rm\, loc.\, sect.}\\ {\rm of}\,\Theta}}} \tilde{{\cal E}}_{\rho(\theta)} /_{(\varphi_\theta^{\theta'})}\textmap{\varphi_{\Theta,\rho}} \mathcal{O}_{U_\theta\times X}\bigg)\textmap{\simeq}(\tilde{\cal E}\textmap{\varphi}{\cal O}_{S\times X}).
$$
\end{pr}

\subsubsection{The natural isomorphism \texorpdfstring{$\id\to \tensor[_S]{\Tg}{}\circ\tensor[_S]{\Ag}{}$}{tgsags-sub}}

 Let $(\Theta,\rho)$ be an object of the groupoid $\Torsss_{G_D}^{\alpha_D}$. We construct a {\it canonical} isomorphism 
$$
(\Theta,\rho)\textmap{\simeq\eta}  \tensor[_S]{\Tg}{}\circ\tensor[_S]{\Ag}{}(\Theta,\rho).
$$ 

In other words, putting $\tilde{\mathcal{E}}\edf \tilde{\mathcal{E}}_{(\Theta,\rho)}$, $\varphi\edf\varphi_{(\Theta,\rho)}$, we will construct a {\it canonical} isomorphism of $\mathcal{O}(G_D)$-torsors $\eta: \Theta\to \Theta_\varphi^D$ making commutative the triangle
$$
\begin{tikzcd}[row sep=5mm]
\Theta\ar[dr, "\rho"']\ar[rr, "\eta \simeq"] & &\Theta^D_\varphi \ar[dl, "\rho^D_\varphi"]\\
&{\cal O}_S(V_D)&
\end{tikzcd}.
$$
For any open subspace $U\subset S$ and $\theta\in \Theta(U)$, denote by 
$$j_2^\theta:{\cal O}_{U\times X} \to p_X^{U\times X*}({\cal E})(D_U)\oplus {\cal O}_{U\times X}$$
 the obvious monomorphism. Using the isomorphisms introduced at the end of section \ref{functorAGS},  we put  
$$
\eta(\theta)\edf (\Psi_\theta\circ j_2^\theta, \psi_\theta).
$$
This formula defines an   isomorphism of $\mathcal{O}(G_D)$-torsors. This follows   using the compatibility formulae (\ref{Compat-psi-S}), (\ref{Compat-Psi-S}) and the commutative diagram (\ref{comm-diag-Psi-theta1-S}). The required equality  $\rho_\varphi^D\circ \eta=\rho$ follows    using the commutative diagram (\ref{comm-diag-Psi-theta2-S}). 
\vspace{1mm}

The collections of functors 	$\tensor[_S]{\Ag}{}$, $\tensor[_S]{\Tg}{}$ associated with complex spaces $S$ define functors 
$$\mathscr{A}:[V_D/G_D]\to \mathscr{A}{\rm ff}_E,\ \mathscr{T}: \mathscr{A}{\rm ff}_E\to [V_D/G_D]$$
over the category $\mathscr{A}\hspace{-0.5mm}{\rm nn}$ of complex spaces. These functors commute with pull-backs and define an equivalence of categories. In conclusion, using the terminology of \cite[p. 11]{Gom}, we have:

\begin{thr}\label{equiv-stacks-th}
Let $X$ be a compact complex manifold, and $E$ be a fixed holomorphic vector bundle on $X$ and ${\cal E}$ the associated locally free sheaf. Let $D$ be an effective divisor on $X$ such that $H^1(X,{\cal E}(D))=0$. The  holomorphic stack $\mathscr{A}{\rm ff}_E$ is isomorphic to the quotient stack $[V_D/G_D]$.
\end{thr}

%
%

%

\end{document}